\documentclass[a4paper,fleqn]{cas-sc}

\usepackage[authoryear]{natbib}
\usepackage{stmaryrd}
\usepackage{placeins}
\usepackage{amsmath}

\def\tsc#1{\csdef{#1}{\textsc{\lowercase{#1}}\xspace}}
\tsc{WGM}
\tsc{QE}
\tsc{EP}
\tsc{PMS}
\tsc{BEC}
\tsc{DE}

\begin{document}
\let\WriteBookmarks\relax
\shorttitle{Elastic wave propagation in fractured media with spring-type and frictional contact deformation laws}
\shortauthors{I.K. Jacobsen et~al.}



\title [mode = title]{Elastic wave propagation in fractured media with spring-type and frictional contact deformation laws}                      

\author[1]{Ingrid Kristine Jacobsen}[orcid=0000-0002-2879-6691]
\cormark[1]
\ead{Ingrid.Jacobsen@uib.no}

\author[1,3]{Jan Martin Nordbotten}[orcid=0000-0003-1455-5704]

\author[1]{Ivar Stefansson}[orcid=0000-0001-6370-496X]

\author[3]{Barbara Wohlmuth}[orcid=0000-0001-6908-6015]

\author[1]{Inga Berre}[orcid=0000-0002-0212-7959]

\affiliation[1]{organization={Center for Modeling of Coupled Subsurface Dynamics, Department of Mathematics, University of Bergen},
                addressline={Allégaten 41}, 
                city={Bergen},
                postcode={5020}, 
                country={Norway}}

\affiliation[2]{organization={NORCE Norwegian Research Centre AS},
                addressline={Nygårdsgaten 12}, 
                city={Bergen},
                postcode={5008}, 
                country={Norway}}

\affiliation[3]{organization={TUM School of Computation, Information and Technology, Technical University of Munich},
                addressline={Boltzmannstraße 3}, 
                city={Garching bei München},
                postcode={85748}, 
                country={Germany}}
                
\cortext[cor1]{Corresponding author}

\begin{abstract}
Elastic wave propagation in fractured media is relevant to applications such as analysis of seismic waves and non-destructive characterization of materials. Understanding attenuation and scattering behavior arising from wave-fracture interaction is important for interpreting observations at both field and laboratory scales. 
This work presents a computational method for elastic wave propagation in fractured media based on a mixed-dimensional discrete fracture-matrix representation. Fracture deformation is governed by four models of increasing complexity, ranging from widely used spring-based formulations to frictional fracture contact models. Many previous studies are restricted to simplified wave fields, single fractures or subsets of the relevant fracture deformation mechanisms. In contrast, the proposed methodology enables fully coupled simulation of elastic wave propagation with fracture deformation models that account for elastic normal deformation, frictional contact and fracture opening and closure.
The elastic wave equation is discretized in space using the cell-centered finite volume method Multi-Point Stress Approximation with weak symmetry and in time using the Newmark method. The spatial discretization is locally conservative and applicable to general polyhedral grids, making it well suited for media containing fractures, material heterogeneities and anisotropy. 
The proposed methodology is verified through numerical convergence analyses and is subsequently applied to wave propagation and fracture deformation in two- and three-dimensional media containing multiple intersecting fractures, demonstrating the influence of fracture constitutive behavior on wave transmission, scattering and fracture deformation.
\end{abstract}

\maketitle

\section{Introduction}\label{sec:introduction}
Wave propagation through fractured media has been extensively studied (e.g. \citet{Schoenberg1980, PyrakNolte1990, Zhao2001, Basabe2016, Hu2024}) due to its relevance for applications such as subsurface seismic analysis \citep{Basabe2016} and non-destructive material characterization \citep{Aleshin2018, Delrue2018}. Seismic waves are scattered and attenuated as they propagate through fractured media \citep{Schoenberg1980, PyrakNolte1990}, and there are several approaches to incorporating the effects of fractures in modeling of wave propagation. Equivalent medium theories (e.g. \citet{Hudson1980}) assign material parameters that represent the averaged effect of multiple fractures within a rock and are well suited for modeling fracture-induced anisotropy. Other approaches assign distinct material parameters locally to account for individual fractures (e.g. \citet{Saenger2002, Parastatidis2021}). However, these approaches generally do not explicitly resolve fracture deformation or local fracture-rock interactions, limiting their ability to capture phenomena such as stick-slip behavior and fracture opening and closure. For a more accurate representation of fracture deformation, as well as fracture-rock interactions, fractures can instead be represented explicitly as discontinuities in the domain. Such representations typically treat fractures as internal boundaries or co-dimension-one objects embedded in the surrounding rock.

Fracture deformation is generally decomposed into tangential and normal components. Both may be modeled using elastic spring-type relations, with either a linear or nonlinear response. Certain applications require the inclusion of frictional effects, in which case tangential deformation may be governed by, for example, Coulomb friction. To allow for fracture opening and closure while ensuring that the fracture surfaces do not penetrate each other, a non-penetration condition may be imposed in the normal direction. The most general fracture deformation models incorporate all of these features.

Among the simplest fracture deformation models are linear spring-type formulations. A commonly used model, introduced by \citet{Schoenberg1980}, relates fracture traction linearly to the displacement jump across the fracture in both the normal and tangential directions. This formulation has been widely applied in analytical and numerical studies of seismic wave propagation in fractured media \citep{Gu1996, Zhang2009, Basabe2016, Vamaraju2018, RamosBarreto2025} and is supported by laboratory experiments \citep{Myer1988, PyrakNolte1990}.

A limitation of linear spring-type models is that they are generally valid only for small wave amplitudes, as larger amplitudes may induce nonlinear fracture deformation \citep{Yi1997}. Several experimental studies have demonstrated that fractures deform nonlinearly (e.g., \citep{Bandis1983, Zhao2008Dynamic}), where the nonlinear deformation is commonly described using the hyperbolic Barton--Bandis relation \citep{Bandis1983}. \citet{Zhao2001} incorporated this relation into the normal component of a spring formulation to study the transmission of compressional waves through a single fracture. The approach has later been employed in numerical simulations of compressional waves propagating through one or several parallel fractures \citep{Hu2024, Zhao2008}. However, these studies are restricted to normally incident compressional plane waves and therefore do not capture tangential fracture deformation.

Under sufficiently strong excitation, fractures may open and close, motivating the development of contact-based fracture representations. \citet{Kimoto2015} imposed unilateral contact conditions in the normal direction together with either sticking or frictionless behavior in the tangential direction. Similarly, \citet{Blanloeuil2020} combined unilateral contact with linear or quadratic spring laws to model nonlinear deformation associated with fracture opening and closure (clapping). However, their analysis was restricted to a one-dimensional configuration with a wave propagating normally across a prestressed fracture. 

Frictional effects have also been incorporated into several fracture deformation models. For example, \citet{Li2013} studied wave propagation through a set of parallel joints using a fracture model combining Barton--Bandis spring-type normal deformation with linear elastic tangential deformation, where the tangential response transitioned to Coulomb slip once a critical tangential displacement was reached. The formulation accounted for obliquely incident compressional and shear waves but was restricted to non-curved wave fronts. \citet{Blanloeuil2014} modeled fractures by an interface of unilateral contact with Coulomb friction to investigate coupled normal and tangential fracture behavior under oblique wave incidence, though without elastic deformation during closure or prior to slip. Alternative formulations for fracture contact with friction include approaches based on memory diagrams \citep{Aleshin2018, Delrue2018} and penalty methods \citep{Gao2019}.

Although these studies capture a range of nonlinear fracture responses, they are often limited by simplified wave fields \citep{Li2013}, sequential coupling strategies \citep{Aleshin2018, Delrue2018, Gao2019}, single-fracture configurations \citep{Blanloeuil2014, Delrue2018, Gao2019}, or fracture models that do not simultaneously account for elastic normal deformation, frictional contact and fracture opening and closure (e.g. \citet{Li2013, Blanloeuil2014}). To our knowledge, fully coupled simulation of elastic wave propagation that simultaneously captures nonlinear fracture deformation, frictional contact and fracture opening and closure has not previously been demonstrated.

In this work, we adopt a monolithic approach that accommodates both spring-type and contact-type fracture deformation models with friction. For the spring-type models, linear \citep{Schoenberg1980} and nonlinear \citep{Bandis1983} stiffness relations are used to relate traction to displacement jump. For the contact-type models, the same stiffness relations govern normal deformation while the fracture surfaces remain in contact, whereas tangential deformation is governed by a Coulomb friction model. This yields a family of fracture deformation models with varying physical complexity, ranging from linear spring-type formulations to nonlinear normal and tangential contact formulations with friction. The fracture deformation models are represented by their own governing equations and solved monolithically together with the elastic wave equation in the surrounding rock.

The coupled problem of rock and fracture deformation requires suitable numerical discretizations, and a range of numerical methods have been developed for this purpose. Finite element methods have been widely used for fracture mechanics and wave propagation, including the work of \cite{Aleshin2018} and \citet{Delrue2018} on wave propagation in media with rough fractures, and Nitsche formulations for static and dynamic contact problems involving frictionless fractures \citep{Chouly2015_p1, Chouly2015_p2} and fractures governed by Coulomb friction \citep{Chouly2023}. The virtual element method extends finite element techniques to polytopal meshes and has been applied to frictionless and Coulomb-frictional contact problems \citep{Laaziri2025Fricless, Laaziri2025Fric}. Enriched finite element methods (e.g. \citet{Komijani2017}) have been employed for wave propagation in fractured media with the advantage that fractures are represented independently of the underlying mesh, beneficial for simulations where the path of the fracture is not known a priori. Higher-order approaches include the hybrid high-order method \citep{Chouly2020} and the discrete de Rham method \citep{Droniou2026}, both applied to quasi-static contact problems with Tresca friction. Discontinuous Galerkin methods constitute another widely used class of methods for wave propagation in fractured media, such as e.g. \citet{Basabe2016} and the hybrid discontinuous Galerkin--spectral element method by \citet{Vamaraju2018}, which provides reduced computational cost compared to the full discontinuous Galerkin discretization.

In this work, we employ the finite volume Multi-Point Stress Approximation method with weak symmetry (MPSA-W) \citep{Keilegavlen2017}. The method has previously been applied to frictional contact problems in quasi-static settings (e.g. \cite{Berge2020, Stefansson2021, Keilegavlen2021Hydromech}), but has not yet been applied to dynamic problems with deformable fractures. The MPSA-W method provides a cell-centered alternative to discontinuous Galerkin methods while retaining geometric flexibility and natural enforcement of free-surface boundary conditions. An additional advantage of finite volume methods is their compatibility with subsurface flow models. Due to their local conservation properties and robustness in heterogeneous media, finite volume discretizations are widely used in reservoir simulation software (e.g. \citet{Lie2019, Ferguson2018, Keilegavlen2021}). This facilitates the development of solution approaches that avoid inter-model interpolation and enables consistent coupling between physical processes. Although the present work is restricted to purely mechanical problems, this compatibility is relevant for future extensions to coupled poromechanical and seismic simulations. 

In our previous work, the MPSA-W method was combined with the implicit Newmark time integration scheme, resulting in the MPSA-Newmark method \citep{Jacobsen2025}. We showed that the MPSA-Newmark method can handle material anisotropy and heterogeneity while achieving second-order convergence in displacement and near-second-order convergence in tractions.

While \citet{Jacobsen2025} demonstrated the accuracy and robustness of the method, the fractures were represented as traction-free interfaces, resulting in complete wave reflection. Accordingly, fracture deformation, wave transmission across fractures and frictional contact were not considered. The present work removes these limitations by monolithically coupling the MPSA-Newmark discretization with governing equations for deformable fractures within a mixed-dimensional model. The main contributions are the development of an MPSA-Newmark methodology for elastic wave propagation in deformable fractured media, the monolithic coupling of the elastic wave equation with fracture deformation models of different complexity, numerical verification of the resulting formulation and illustrative simulation examples of wave-fracture interaction in two and three dimensions.

The paper is organized as follows. Section 2 presents the mixed-dimensional geometry representation and the governing equations, Section 3 describes the spatial and temporal discretizations and Section 4 contains verification of the methodology through numerical convergence analyses. Section 5 presents a comparison of the four fracture deformation models and Section 6 presents two- and three-dimensional simulation examples in domains containing multiple intersecting fractures. Finally, the conclusion is presented in Section 7.
\FloatBarrier
\section{Mathematical model}\label{sec:mathematical_model}
In this section we present the mathematical model. Specifically, we cover the mixed-dimensional representation of the geometry, the elastic wave equation and the fracture deformation models we consider in this paper. Finally we summarize initial and boundary conditions to close the mathematical model.

\subsection{Mixed-Dimensional Representation of the Geometry}\label{subsec:mixed_timensional_geometry}
The simulation domain $\Omega$ is modeled as a collection of subdomains of different dimensions, where fractures are represented as co-dimension-one objects. A three-dimensional simulation domain may thus consist of a 3D rock subdomain and 2D fracture subdomains. An analogous subdivision is used in two dimensions, with 2D rock and 1D fractures.

\begin{figure}[pos=htbp]
    \centering
    \includegraphics[width=0.5\linewidth]{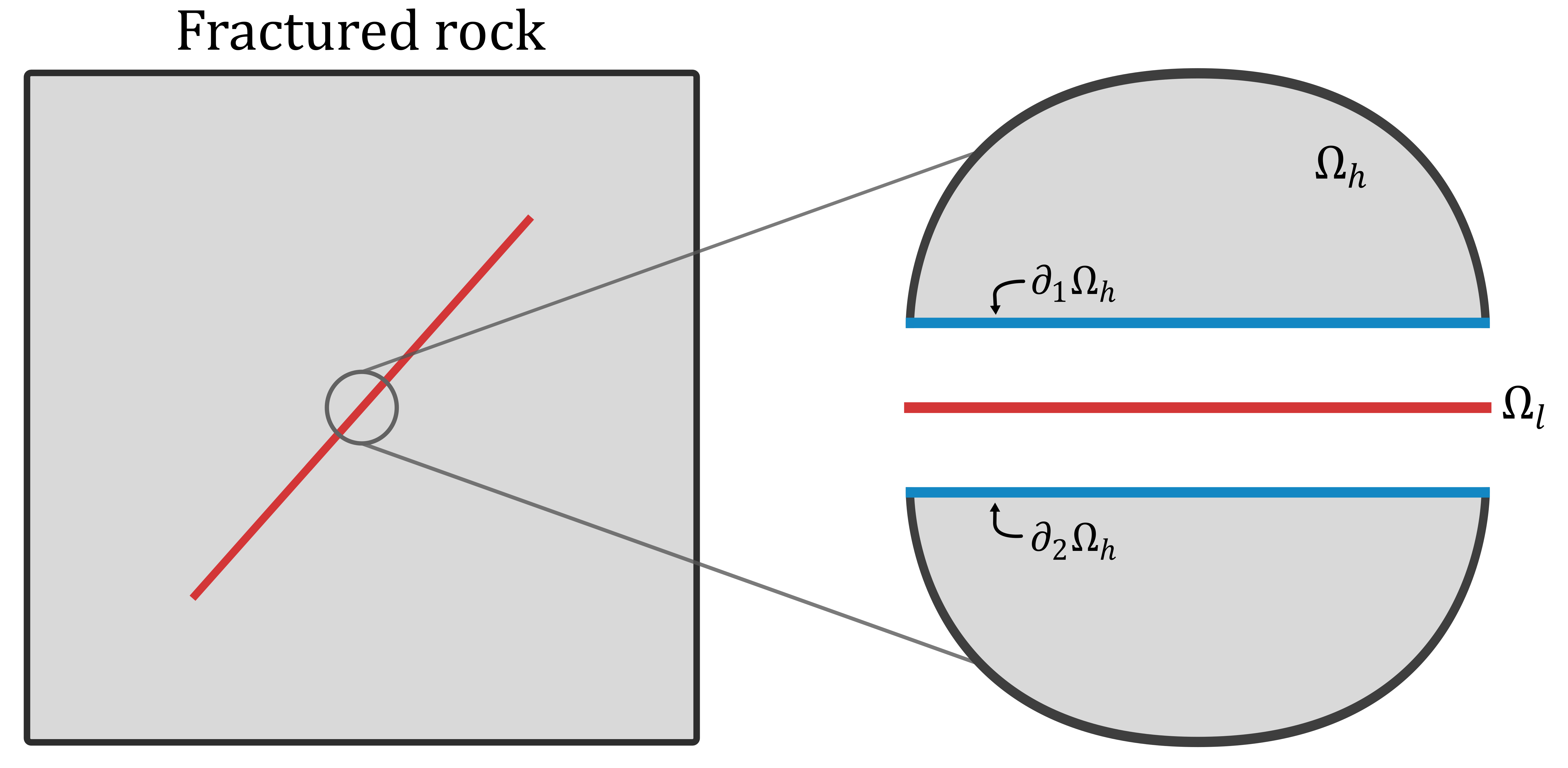}
    \caption{Left: Schematic of a rock ($\Omega_h$, grey) with a fracture ($\Omega_l$, red). Right: A zoom-in showing the fracture and the higher-dimensional internal boundary $\partial_i\Omega_h$ for $i \in \{1,2\}$ (blue), separated for visualization}
    \label{fig:burgermodell}
\end{figure}
Though $\Omega$ is in general 3D, we present the mixed-dimensional representation of the geometry in 2D for illustrative purposes. We therefore refer to Figure~\ref{fig:burgermodell}, which shows a 2D rock $\Omega_h$ with a 1D fracture $\Omega_l$, whose two sides are denoted by $\partial_i\Omega_h$ for $i \in \{1,2\}$. A fracture, such as $\Omega_l$ in Figure~\ref{fig:burgermodell}, is associated with a unit normal vector denoted by $n_l$. By convention, this vector is aligned with the outward unit normal on side $1$ of the fracture and points in the opposite direction of that on side $2$. Consequently, the dot product between $n_l$ and the outward unit normal on side $i$ satisfies: $n_l\cdot n_i=(-1)^{i-1}$ for $i\in\{1, 2\}$. 

The jump of a quantity across a fracture is defined as
\begin{equation}
    \llbracket \cdot \rrbracket = (\cdot)_2 - (\cdot)_1,
        \label{eq:jump}
\end{equation}
where the subscripts $1$ and $2$ denote quantities on $\partial_1 \Omega_h$ and $\partial_2 \Omega_h$, respectively. These will also be referred to as interface quantities (e.g. interface displacements). We adopt the convention that compressive stresses are negative. We likewise define compressive normal displacement jumps to be negative, so that they correspond to fracture closure.

Finally, we introduce the notion of normal and tangential components of a generic vector $\iota$ associated with a fracture:
\begin{equation}
    \iota_n=\iota \cdot n_l, \qquad \iota_\tau=\iota-\iota_nn_l.
\end{equation}
Note that $\iota_\tau$ denotes a vector with either one or two components, depending on whether $\iota$ is a two- or three-dimensional vector, respectively.

\subsection{Governing Equations in the Rock Matrix}\label{subsec:governing_equations_rock_matrix}
Conservation of momentum in the rock matrix is represented by the elastic wave equation:
\begin{equation}
    \rho \ddot u = \nabla \cdot \sigma + h.
    \label{eq:wave_eq}
\end{equation}
Here $\rho$ is the rock density, $h$ is an external force term and, using the dot notation for time-derivatives, $\ddot u$ denotes the acceleration. $\sigma$ denotes the stress tensor which is governed by Hooke’s law:
\begin{equation}
    \sigma = \mathcal C : \epsilon\left(u\right).
    \label{eq:hookes}
\end{equation}
Here $\epsilon$ denotes linearized strain:
\begin{equation}
    \epsilon = \frac{1}{2} \left( \nabla u + (\nabla u)^T \right),
    \label{eq:symmetric_gradient}
\end{equation}
and $\mathcal{C}$ denotes the fourth order stiffness tensor. Using tensor notation, $\mathcal{C}$ for an isotropic solid is formulated as
\begin{equation}
    \mathcal{C}_{ijkl}=\lambda \delta_{ij} \delta_{kl} + \mu(\delta_{ik}\delta_{jl} + \delta_{il}\delta_{jk}).
\end{equation}
where $\lambda$ and $\mu$ denote the first and second Lamé parameter, respectively, and $\delta_{ij}$ is the Kronecker delta taking the value 1 if $i=j$ and 0 otherwise.

\subsection{Fracture Deformation Models}\label{subsec:fracture_deformation_models}
Two model types are considered: the spring model (S) and the contact mechanics model (C). In both cases, fracture deformation is decomposed into tangential and normal components. 

In the spring model (S), both the normal and tangential behavior are governed by spring-type deformation relations. The tangential behavior is linear (Lin), while the normal behavior is either linear (Lin) or described by the nonlinear Barton--Bandis (BB) relation. 

In the contact mechanics model (C), tangential deformation is governed by Coulomb friction (Coul), while a unilateral non-penetration constraint is enforced in the normal direction. Normal elastic deformation is activated only when contact occurs and follows either a linear (Lin) or Barton--Bandis (BB) relation, consistent with the spring model. 

The model variants are denoted according to the naming convention ``model type - tangential deformation - normal elastic deformation'', resulting in the four combinations S-Lin-Lin, S-Lin-BB, C-Coul-Lin and C-Coul-BB:
\begin{itemize}
    \item S-Lin-Lin: Linear deformation in both tangential and normal directions.
    \item S-Lin-BB: Linear tangential deformation with nonlinear Barton–Bandis normal deformation.
    \item C-Coul-Lin: Contact mechanics model with Coulomb friction in the tangential direction and linear normal deformation during contact, subject to a non-penetration constraint.
    \item C-Coul-BB: Contact mechanics model with Coulomb friction in the tangential direction and nonlinear Barton–Bandis normal deformation during contact, subject to a non-penetration constraint.
\end{itemize}
The details of the individual fracture deformation models are presented below.

\subsubsection{Fracture contact traction}\label{subsubsec:fracture_contact_traction}
The fracture deformation models are formulated in terms of fracture contact traction $q_l$ which, by convention, is taken as the surface traction on side $1$ of the fracture. With traction defined as $q=\sigma \cdot n$, the surface traction on fracture side $i$ is $q_i=\sigma_i\cdot n_i$ for $i\in\{1,2\}$. Applying Newton's third law yields:
\begin{equation}
    q_l=q_1=-q_2.
    \label{eq:newtons_3rd}
\end{equation}
For notational simplicity, the subscript is omitted and fracture traction is denoted by $q=q_l$ throughout this section.

\subsubsection{Spring-type fracture deformation}\label{subsubsec:spring_type_fracture_deformation}
The spring model assumes constant traction across a fracture, but allows discontinuity in the displacement field. The tangential and normal fracture deformation equations are expressed by: 
\begin{align}
    \llbracket u \rrbracket_\tau - g_\tau^{\text{Lin}}(q_\tau) &= 0, \label{eq:S_tangential}\\
    \llbracket u \rrbracket_n - g_n^{\eta}(q_n) &= 0 \qquad \text{for} \quad \eta\in\{\text{Lin, BB}\}, \label{eq:S_normal}
\end{align}
where $g_\tau^{\text{Lin}}(q)$ and $g_n^\eta(q)$ denote the tangential and normal stiffness relations, respectively. The tangential relation is always linear, whereas the normal relation may be either linear or nonlinear. The specific forms of $g_\tau^{\text{Lin}}(q)$ and $g_n^\eta(q)$ are given in Section~\ref{subsubsec:stiffness_relations}.

\subsubsection{Fracture contact mechanics with friction}\label{subsubsec:fracture_contact_mechanics}
The second fracture deformation model considered is contact mechanics with friction. Coulomb friction is assumed in the tangential direction, and the friction bound is denoted by
\begin{equation}
    b=-Fq_n,
\end{equation}
where $F$ is the coefficient of friction. The friction model is given by
\begin{equation}
    \begin{aligned}
        \lVert q_\tau \rVert &\leq b, \\
        \lVert q_\tau \rVert &< b \rightarrow \dot {\llbracket u \rrbracket}_\tau = 0, \\
        \lVert q_\tau \rVert &= b \rightarrow \dot {\llbracket u \rrbracket}_\tau = k q_\tau, \quad k > 0,
    \end{aligned}
    \label{eq:friction_law}
\end{equation}
where $\dot {\llbracket u \rrbracket}_\tau$ is the time derivative of the tangential displacement jump, hereafter referred to as the tangential velocity jump. For future reference we define $s=\lVert q_\tau \rVert/\lvert b\rvert$ as the slip tendency of a fracture: $s=1$ means slip whereas $s<1$ means no slip.

In the normal direction, the following non-penetration condition is imposed:
\begin{equation}
\begin{aligned}
    \llbracket u \rrbracket_n - g_n^{\eta}(q_n) &\geq 0 \\
    q_n\left({\llbracket u \rrbracket}_n-g_n^{\eta}(q_n) \right) &= 0 \qquad \qquad \text{for}\quad \eta\in\{\text{Lin, BB}\}.\\
    q_n &\leq 0
    \label{eq:non_penetration}
\end{aligned}
\end{equation}
In this context, $g_n^{\eta}(q_n)$ represents the elastic normal deformation of the fracture when the fracture surfaces are in mechanical contact. For future reference, the fracture opening is defined as $\delta = \llbracket u \rrbracket_n - g_n^{\eta}(q_n)$.

The contact mechanics models considered herein yield one of the following contact states at each location on the fracture:
\begin{itemize}
    \item Open: The fracture sides are not in mechanical contact ($\delta>0$, $q=0$, $s = \mathrm{NaN}$).
    \item Sticking: The fracture sides are in mechanical contact and not moving relative to each other ($q\neq0$, $\lVert q_\tau\rVert<b$, $s<1$). 
    \item Sliding: The fracture sides are in mechanical contact and are moving relative to each other ($q\neq0$, $\lVert q_\tau\rVert=\lvert b \rvert$, $s=1$).
\end{itemize}

Enforcing the contact mechanics conditions is done by reformulating the Inequalities~\ref{eq:friction_law}~and~\ref{eq:non_penetration} into equalities using the so-called radial-return map formulation \citet{Alart1991}:
\begin{align}
    q_\tau - \min\left\{1, \frac{\max\{b, 0\}}{\lVert q_\tau + c \dot{\llbracket u \rrbracket}_\tau \rVert}\right\}
    (q_\tau + c \dot{\llbracket u \rrbracket}_\tau)
    &= 0,
    \label{eq:C_tangential}
    \\[1ex]
    q_n - \min\left\{0, q_n + c\left({\llbracket u \rrbracket}_n - g_n^{\eta}(q_n)\right)\right\}
    &= 0
    \quad \text{for } \eta\in\{\text{Lin, BB}\},
    \label{eq:C_normal}
\end{align}
where $c>0$ is a numerical parameter.

\subsubsection{Stiffness relations}\label{subsubsec:stiffness_relations}
The stiffness relations define the relationship between displacement jump and fracture contact traction in Equations~\ref{eq:S_tangential}, \ref{eq:S_normal} and \ref{eq:C_normal}. In the normal direction, either a linear \citet{Schoenberg1980} or a nonlinear \citet{Bandis1983} relation is employed:
\begin{equation}
    g_n^{\text{Lin}}(q_n) = \frac{q_n}{K_n}, \qquad g_n^{\text{BB}}(q_n) = \frac{q_n}{K_n-\frac{q_n}{\Delta u_{\max}}},
    \label{eq:g_BB}
\end{equation}
whereas in the tangential direction, we consider only the linear relation:
\begin{equation}
    g_\tau^{\text{Lin}}(q_\tau) = \frac{q_\tau}{K_\tau}.
\end{equation}
Here, $K_n$ and $K_\tau$ denote the normal and tangential fracture stiffnesses, respectively, and $\Delta u_{\max}$ is the maximum allowed fracture closure. We emphasize that $g_\tau^{\mathrm{Lin}}$ does not enter the contact mechanics model, where tangential deformation is governed by the Coulomb friction law. Hence, $g_\tau^{\mathrm{Lin}}$ appears only in the spring models.

\subsection{Boundary and Initial Conditions}\label{subsec:bounday_and_initial_conditions}
Initial and boundary conditions are specified to close the system. The outer boundary of the domain is denoted by $\Gamma$ and partitioned into three non-overlapping parts, such that $\Gamma=\Gamma_D\cup\Gamma_N\cup\Gamma_A$. These correspond to Dirichlet (D), Neumann (N) and absorbing (A) boundaries.

Neumann conditions are set in terms of traction boundary values $F_N$:
\begin{equation}
    \sigma \cdot n = F_N \quad \text{on } \Gamma_N,
\end{equation}
where $n$ is the outward pointing normal vector relative to the boundary. Dirichlet conditions are set in terms of displacement boundary values $F_D$:
\begin{equation}
    u = F_D \quad \text{on } \Gamma_D.
\end{equation}
The absorbing boundary conditions \citep{Clayton1977, Tsogka1999}, following the presentation in \citet{Jacobsen2025}, can be expressed as:
\begin{equation}
    \sigma \cdot n + \mathcal{D}\,\dot{u} = 0 \quad \text{on } \Gamma_A,
\end{equation}
where $\dot u$ denotes the velocity and $\mathcal D = \left(\sqrt{\rho(\lambda + 2\mu)}nn^T + \sqrt{\rho\mu}(I-nn^T)\right)$.
\FloatBarrier
\section{Discretization}\label{sec:discretization}
This section describes the spatial and temporal discretization of the model presented in the previous section, as well as the resulting system to be solved. The model is discretized in space by the finite volume Multi-Point Stress Approximation with weak symmetry and in time by the Newmark method. 

\subsection{Multi-Point Stress Approximation with Weak Symmetry}\label{sec:mpsa_w}
Several presentations of the finite volume Multi-Point Stress Approximation with weak symmetry (MPSA-W) are available in the literature, and we therefore limit the discussion herein, referring to \citet{Keilegavlen2017}, \citet{Nordbotten2021} and \citet{Berge2020} for details. In particular, Section 3 (excluding Section 3.2.1) of \citet{Berge2020} provides an in-depth description of MPSA-W for fractured media. We note that \citet{Berge2020} considers poroelastic media, and only the parts concerning mechanical deformation are relevant to the present work.

\begin{figure}[pos = htbp]
    \centering
    \includegraphics[width=0.3\linewidth]{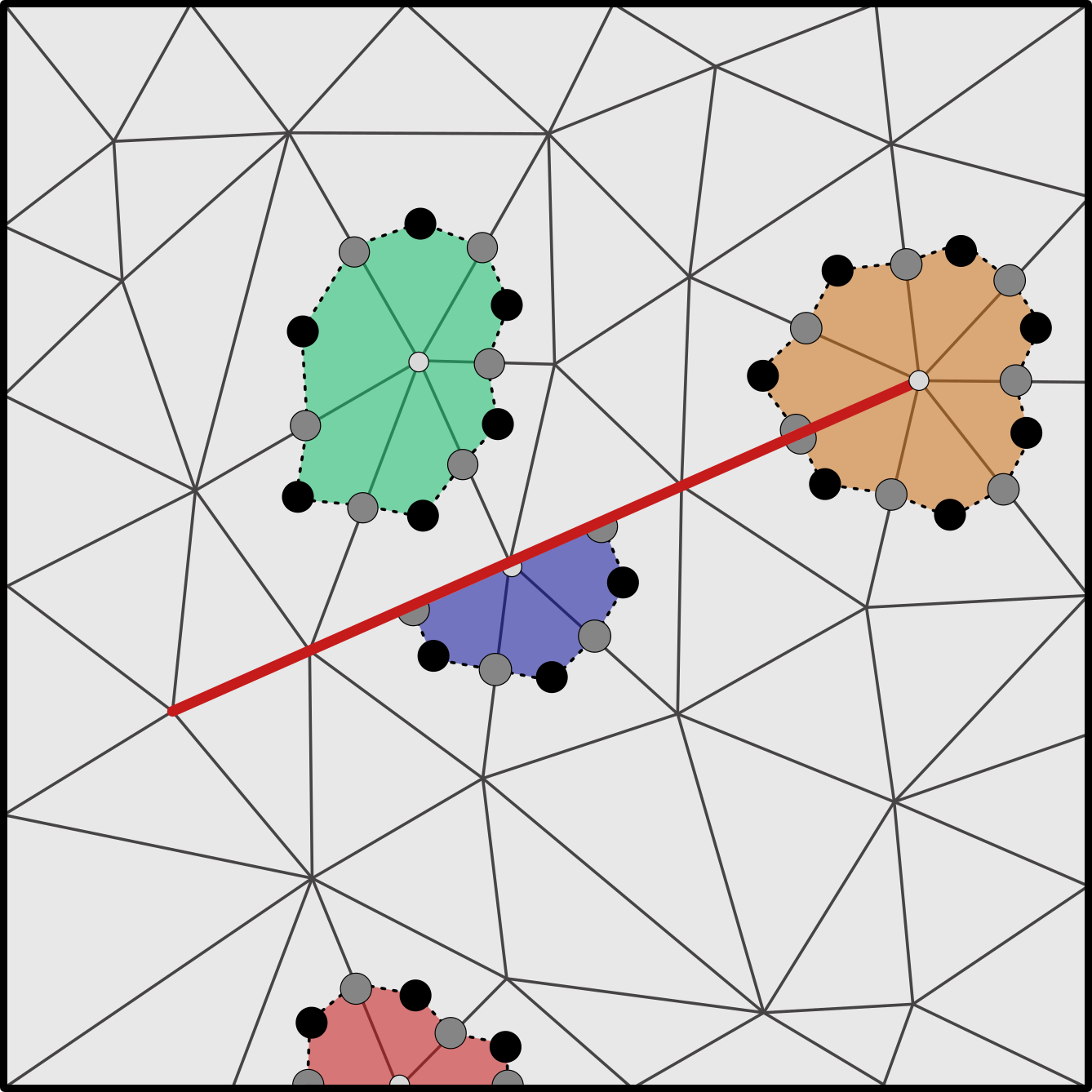}
    \caption{Grid with a diagonal fracture (red line). The colored patches show different interaction regions for the discretization. Cell-centers, face-centers and cell-vertices that are involved in the interaction regions are denoted by black, dark gray and light gray circles, respectively. Half circles indicate that the interaction region involves a boundary or a fracture}
    \label{fig:mpsa_interaction_regions}
\end{figure}

\begin{figure}[pos = htbp]
    \centering
    \includegraphics[width=0.4\linewidth]{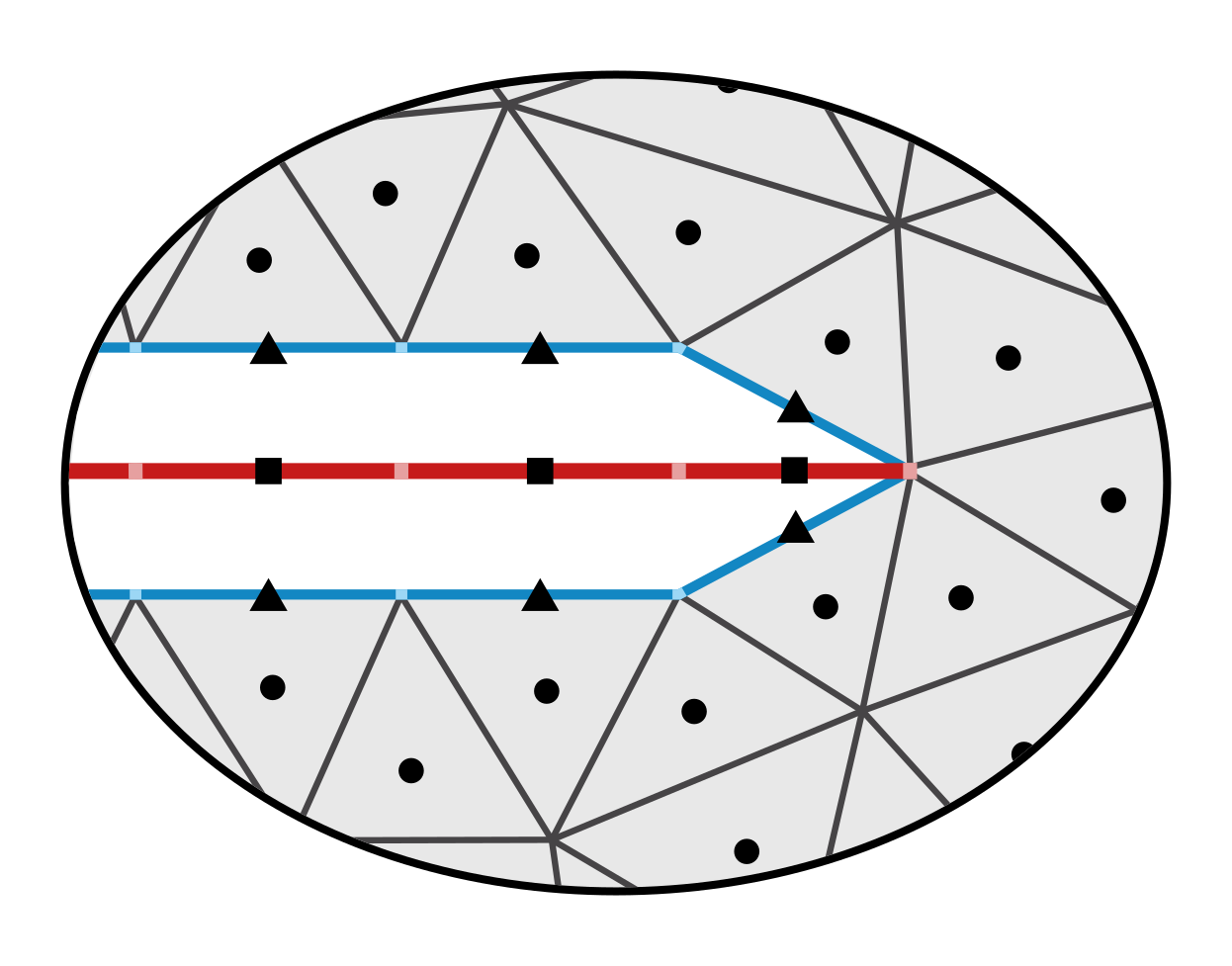}
    \caption{Schematic of the grid for the fracture and matrix and the locations of the discrete unknowns. The rock matrix grid is shown as grey triangles, while the fracture grid is represented by a red horizontal line with pink divisions indicating the faces between fracture cells. The blue lines denote the internal boundary of the rock matrix grid, with light blue indicating a partition which is aligned with the fracture cells and faces of the higher dimensional grid. Black circles: cell-centered displacements $\boldsymbol{u}$. Black triangles: interface displacements $\boldsymbol{u}_1$ and $\boldsymbol{u}_2$. Black squares: fracture contact tractions $\boldsymbol{q}_l$}
    \label{fig:dof_location}
\end{figure}

The MPSA-W method is based on continuity conditions for displacement and stress, which are used to construct local linear systems on interaction regions around each grid vertex. Solving these systems gives a mapping to displacement gradients expressed from cell-center displacements and, where applicable, boundary values and fracture unknowns. Boundary values enter the local systems when the interaction region intersects the boundary, whereas if the interaction region intersects a fracture, the associated fracture contact tractions enter the local systems as additional unknowns. We refer to Figure~\ref{fig:mpsa_interaction_regions} for an illustration of the fractured domain and the four types of interaction regions, and to Figure~\ref{fig:dof_location} for a schematic showing the location of the discrete unknowns. After the displacement gradients are obtained, they can be eliminated from the system, resulting in a system with cell-center displacements and fracture contact tractions as the remaining unknowns.

\subsection{The Newmark Method}\label{sec:newmark_method}
The model is discretized in time using the Newmark method \citep{Newmak1959}, which in its original formulation expresses the unknown displacement and velocity as functions of the unknown acceleration. To remain consistent with the primary unknowns in the spatial discretization, namely displacements, we do not introduce velocity or acceleration as new unknowns to the model. We therefore use the Newmark method in its displacement formulation, where velocity and acceleration are expressed as functions of unknown displacements \citep{Chopra2012, Newmak1959}
\begin{align}
    {\dot u}^n &= \left(1 - \frac{\gamma}{\beta}\right){\dot u}^{n-1} + \Delta t \left(1 - \frac{\gamma}{2\beta}\right) {\ddot u}^{n-1} + \frac{\gamma}{\beta \Delta t}(u^n - u^{n-1}), 
    \label{eq:newmark_velocity}\\
    {\ddot u}^n &= \frac{1}{\beta \Delta t^2}\left(u^n - u^{n-1} - \Delta t {\dot u}^{n-1} - (1 - 2\beta)\frac{\Delta t^2}{2}{\ddot u}^{n-1}\right).
    \label{eq:newmark_acceleration}
\end{align}
Here the superscripts $n$ and $n-1$ denote the current and previous time step number, respectively, and $\beta$ and $\gamma$ are discretization parameters. We set $\beta=1/4$ and $\gamma=1/2$ which is a standard choice (see e.g. \citet{Subbaraj1989}). 

The time discretization is done by substituting the Newmark expressions into the time derivatives in the model equations. Specifically, substituting Equation~\ref{eq:newmark_acceleration} into Equation~\ref{eq:wave_eq} leaves the current displacement as the only unknown arising from the acceleration term and with the known quantities $u^{n-1}$, $\dot u^{n-1}$ and $\ddot u^{n-1}$ entering as source-like terms. A similar procedure is done for the velocity jump in Equation~\ref{eq:C_tangential}. Due to linearity of the jump operator, Equations~\ref{eq:newmark_velocity}~and~\ref{eq:newmark_acceleration} can also be applied to jumps in quantities across fractures. Equation~\ref{eq:newmark_velocity} is therefore used to eliminate the velocity jump in Equation~\ref{eq:C_tangential}. The remaining unknown from the time-discretized velocity jump is thus the jump in interface displacements.

At the end of each successful time step, when the unknown displacements are obtained, the velocity and acceleration, as well as their corresponding jumps, are updated according to Equations~\ref{eq:newmark_velocity}~and~\ref{eq:newmark_acceleration}.

\subsection{Fully Discrete System and Solution Procedure}\label{sec:fully_discrete_system}
The nonlinear system to be solved is constructed by the spatially discretized governing equations using MPSA-W (Section~\ref{sec:mpsa_w}), with the temporal discretization applied as described in Section~\ref{sec:newmark_method}. The system consists of the discretized elastic wave equation (Equation \ref{eq:wave_eq}), fracture deformation equations (Equations~\ref{eq:S_tangential}~and~\ref{eq:S_normal} or Equations~\ref{eq:C_tangential}~and~\ref{eq:C_normal}), continuity of tractions (Equation~\ref{eq:newtons_3rd}) and boundary conditions. With reference to Figure~\ref{fig:dof_location}, the unknowns in the linear system are then cell-centered displacements ($\boldsymbol{u}$, circles), interface displacements ($\boldsymbol{u}_1$ and $\boldsymbol{u}_2$, triangles) and fracture contact tractions ($\boldsymbol{q}_l$, squares).

Solving the system is done by using the Newton-Raphson method. When fracture contact mechanics with friction is considered, that is, when fracture deformation is represented by C-Coul-Lin or C-Coul-BB, the system is only semi-smooth and is therefore solved using a semi-smooth Newton method \citep{Berge2020, Hueber2008}.

All simulations are carried out using the open-source code PorePy \citep{Keilegavlen2021}. The simulations employ iterative linear solvers using GMRES with a restart parameter of 100. The linearized coupled systems are solved using a Schur complement reduction of the fracture equations \citep{Zabegaev2026}, with algebraic multigrid (AMG) applied to the coupled system of elasticity and force balance equations on the fracture surfaces. The AMG preconditioner is based on the HYPRE BoomerAMG implementation, using unknown-based coarsening and strong threshold values of 0.25 for two-dimensional simulations and 0.7 for three-dimensional simulations.

The source code and run scripts for all simulations are available in \citet{Jacobsen2026Zenodo} and were developed by the authors with subsequent AI assisted refactoring and improvements. Plotting scripts for visualization of data were generated using generative AI and subsequently tailored by the authors. All code generated or modified with generative AI assistance was reviewed and verified by the authors. The generative AI tool used was ChatGPT.

\FloatBarrier
\section{Numerical Convergence Analyses}\label{sec:numerical_convergence_analyses}
This section presents numerical convergence analyses of wave propagation through fractures modeled using the fracture deformation models introduced in Section~\ref{sec:mathematical_model}. Unless otherwise specified, all simulations in this section are performed using the domain geometry and material parameters described below. In addition, all discretization grids are unstructured triangular grids; see Figure~\ref{fig:domain_A_B} for an example.

The domain and material parameters are based on those used in \citet{Gao2019}, except for the fracture stiffnesses, which were not part of their model. The domain is rectangular, measuring 50\,mm by 25\,mm, and the material represents a soft rock with density $\rho = 2600\,\mathrm{kg\,m^{-3}}$ and Lamé parameters $\mu=\lambda=4.0\,\mathrm{GPa}$. We choose  fracture stiffnesses $K_n=K_\tau=2.0\cdot10^{11}\,\mathrm{Pa\,m^{-1}}$ for all simulations in this section. The boundary conditions are specified separately for each simulation case.

\subsection{Spring-Type Fracture Deformation: Linear Model}\label{subsec:spring_type_linear}
This section considers elastic waves propagating with normal incidence towards a fracture governed by the S-Lin-Lin model. An analytical solution based on \citet{Schoenberg1980} is used to compute the error in the numerical solution following the procedure described in Appendix~\ref{appendix_known_analytical_solution}.

\subsubsection{Analytical solution}\label{subsubsec:analytical_solution}
\begin{figure}[pos=htbp]
    \centering
    \includegraphics[width=0.4\linewidth]{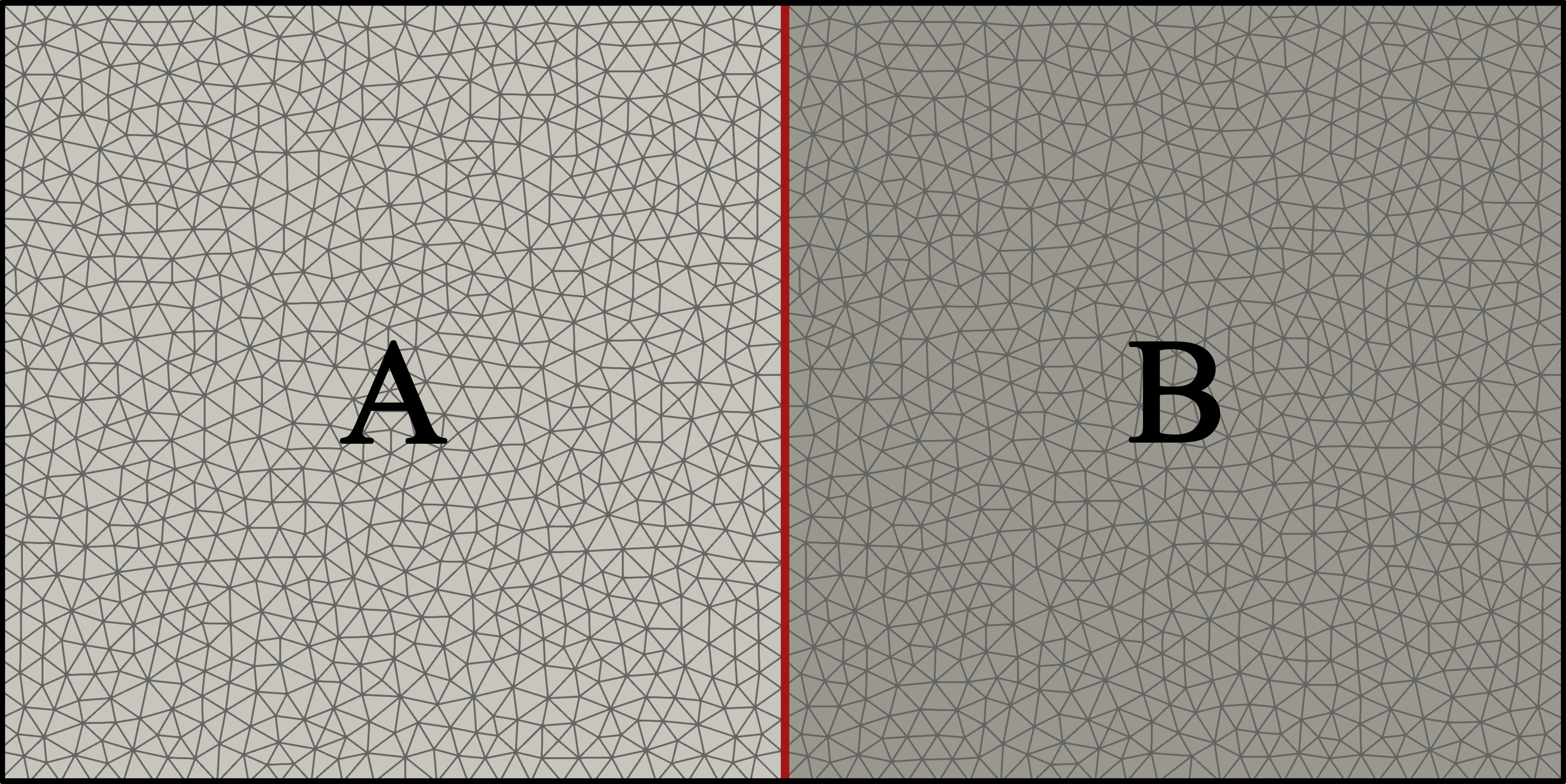}
    \caption{Schematic of the simulation domain used in Sections~\ref{subsec:spring_type_linear} and~\ref{subsec:spring_type_BB}. The domain is divided into regions A and B with potentially different material parameters (e.g. a softer and a stiffer rock). The fracture at the A-B interface is shown as a vertical red line, while the dark gray grid lines exemplify the appearance of the unstructured triangular mesh used in these sections}
    \label{fig:domain_A_B}
\end{figure}

We consider compressive and shear waves propagating at normal incidence towards a vertical fracture located at $x=l$, at the interface between the two dissimilar media A and B (see Figure~\ref{fig:domain_A_B}). The convergence analyses are performed in a 2D domain, although wave propagation occurs only in the $x$-direction.

\begin{figure}[pos=htbp]
    \centering
    \includegraphics[width=0.45\linewidth]{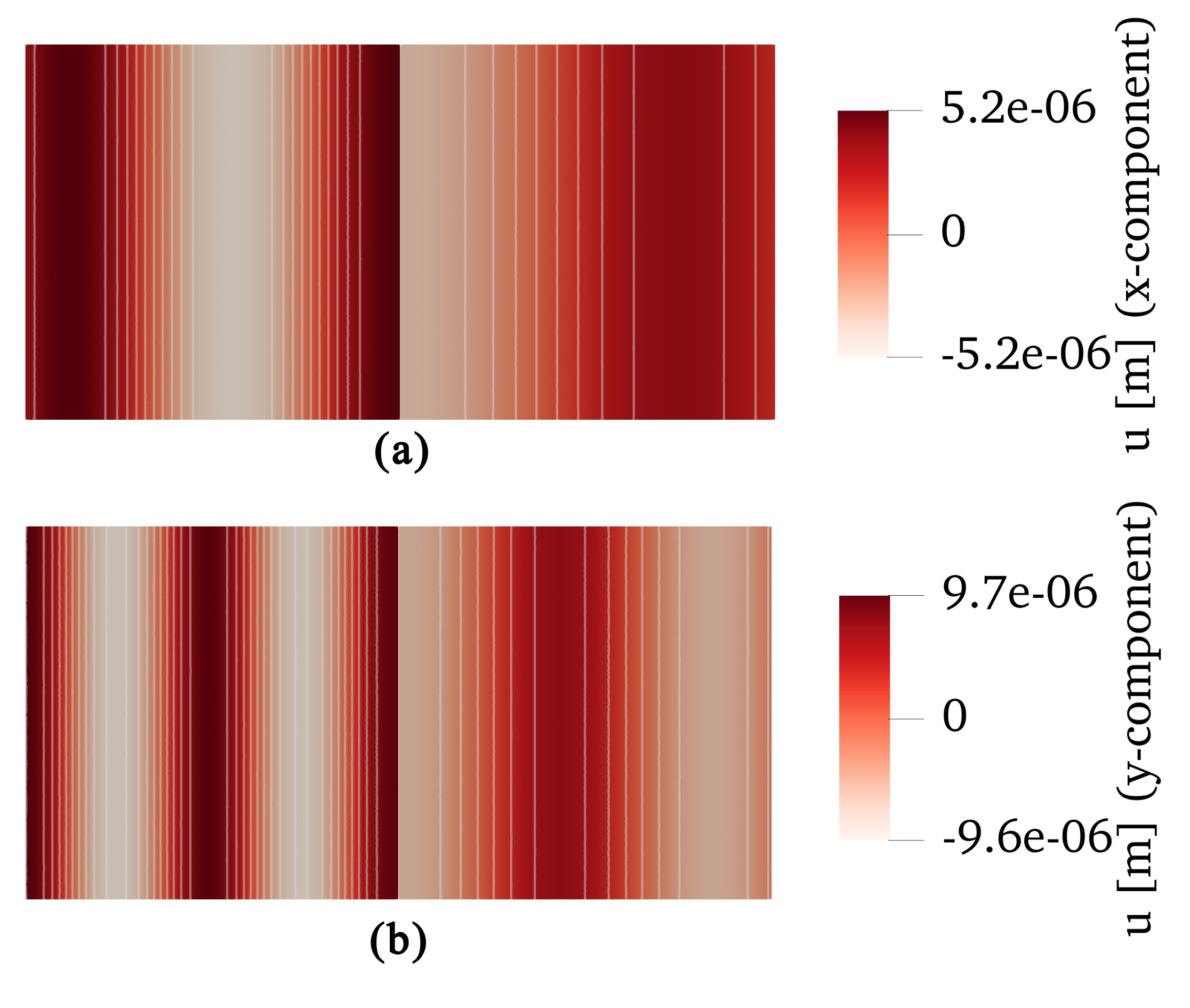}
    \caption{Displacement fields for the convergence analyses with fractures governed by the S-Lin-Lin model. (a) Shows the $x$-component in Equation~\ref{eq:u_p} for the compressive setup. (b) Shows the $y$-component in Equation~\ref{eq:u_s} for the shear setup. In both (a) and (b), a clear discontinuity in the displacement field is observed across the fracture. Contour lines are added for visualization purposes}
    \label{fig:spring_convergence_compressive_and_shear}
\end{figure}

The analytical solutions for the compressive and shear setups can be formulated as 
\begin{align}
    u_p(x, t)&=[u(x, t),0]^T, \label{eq:u_p} \\
    u_s(x, t)&=[0, u(x, t)]^T, \label{eq:u_s}
\end{align}
where subscripts $p$ and $s$ denote compressive and shear, respectively, and the function $u(x, t)$ represents the wave field. We emphasize that although the analytical solutions are one-dimensional, the problem is solved on an unstructured triangular grid (Figure~\ref{fig:domain_A_B}), providing a fully two-dimensional verification of the numerical method. Using standard techniques to obtain the wave field $u(x,t)$ we get:
\begin{align}
    u(x,t) =
    \begin{cases}
        U \sin\!\left(\omega\left(t - \frac{x - l}{c_{A}}\right)\right)
        + U \lvert R(\omega) \rvert  \sin\!\left(\omega\left(t + \frac{x - l}{c_{A}}\right) + \phi_R\right) & \text{for}\, x < l, \\
        U \lvert T(\omega)\rvert  \sin\!\left(\omega\left(t - \frac{x - l}{c_{B}}\right) + \phi_T\right) & \text{for}\, x \ge l,
    \end{cases}
\label{eq:anal_sol_linear_spring}
\end{align}
where $\lvert \cdot \rvert$ and $\phi_{(\cdot)}$ denote the modulus and argument, respectively, of the reflection coefficient $R(\omega)$ and the transmission coefficient $T(\omega)$. The symbols $c_A$ and $c_B$ denote the wave speeds in media A and B, respectively. Whether these correspond to shear or compressive wave speeds is specified below. Figure~\ref{fig:spring_convergence_compressive_and_shear} illustrates the wave fields in the compressive and shear setups.

Following \citet{Schoenberg1980}, the reflection and transmission coefficients $R(\omega)$ and $T(\omega)$ for plane waves interacting with a fracture at the interface between two media A and B are given by
\begin{equation}
    \begin{aligned}
        T(\omega)&=\frac{2Y_{B,(i,j)}}{i+Y_{B,(i,j)}+Y_{A,(i,j)}}\\[1ex]
        R(\omega)&= \frac{i+Y_{B,(i,j)}-Y_{A,(i,j)}}{i+Y_{B,(i,j)}+Y_{A,(i,j)}}
    \end{aligned}
\qquad
\text{for } (i,j)\in\{(n, p), (\tau, s)\},
\end{equation}
with $Y_{m,(i,j)} = K_i/(\omega \rho_m c_{j,m})$ for $m \in \{\text{A},\text{B}\}$. The wave speed $c_{j,m}$ depends on the wave mode. For compressive waves, $c_{p,m} = \sqrt{(\lambda_m + 2\mu_m)/\rho_m}$, whereas for shear waves, $c_{s,m} = \sqrt{\mu_m/\rho_m}$. The parameter $K_i$ represents either the normal or tangential fracture stiffness, depending on the wave type: $K_i=K_n$ for the compressive setup and $K_i=K_\tau$ for the shear setup.

Computing the modulus and argument of $R(\omega)$ and $T(\omega)$ and inserting Equation~\ref{eq:anal_sol_linear_spring} into Equations~\ref{eq:u_p} and \ref{eq:u_s} yields an analytical solution for the propagation of compressive and shear waves, respectively, in a heterogeneous fractured medium.
\subsubsection{Numerical convergence analyses setup and results}\label{subsubsec:numerical_convergence_analyses_linear_spring}
We consider a rectangular domain $\Omega$ with a fracture located at $x=l=25.0\,\mathrm{mm}$. Different material parameters are assigned to media A and B, rendering the domain heterogeneous. Medium A has the material parameters specified at the beginning of this section ($\mu=\lambda=4.0\,\mathrm{GPa}$ and $\rho=2600.0\,\mathrm{kg\,m^{-3}}$). Medium B represents a shale-like material with $\lambda=11.0\,\mathrm{GPa}$, $\mu=10.2\,\mathrm{GPa}$ (converted from the Young's modulus and Poisson ratio reported by \citet{Sun2022}) and a representative density of $\rho=2400.0\,\mathrm{kg\,m^{-3}}$.

\begin{figure}[pos=htbp]
    \centering
    \includegraphics[width=0.99\textwidth]{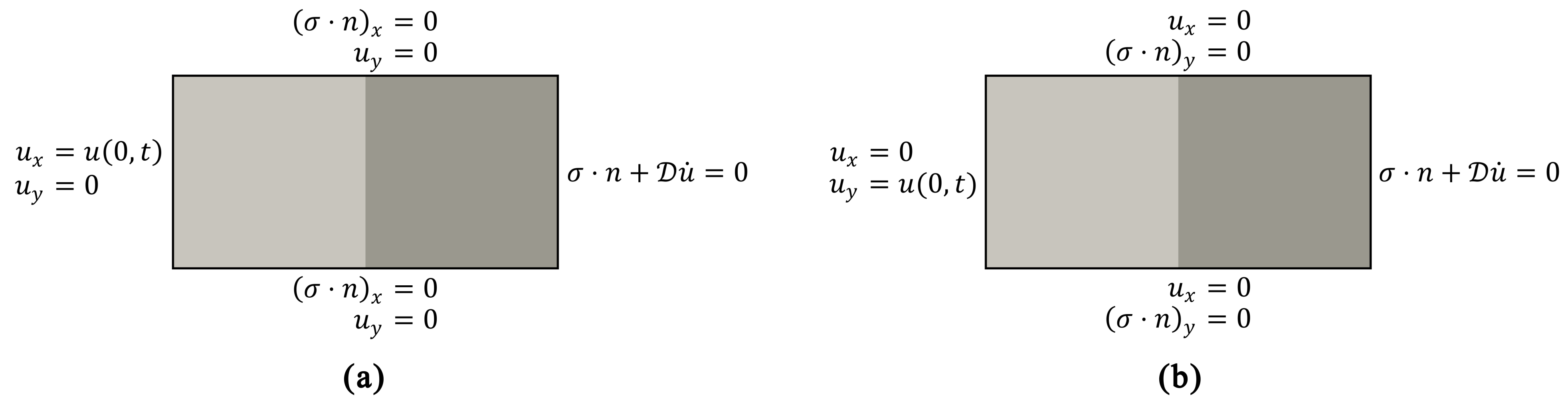}
    \caption{Boundary conditions for the convergence analyses with spring fracture deformation models. The left boundary is subject to a time-dependent Dirichlet condition, the right boundary is absorbing and the top and bottom boundaries are subject to roller conditions. (a) Compressive setup. (b) Shear setup}
    \label{fig:bc_type_comp_shear}
\end{figure}

Initial and boundary conditions are specified to close the system. The initial displacement is obtained from Equation~\ref{eq:anal_sol_linear_spring} evaluated at $t=0$, while its first and second time derivatives are used to define the initial velocity and acceleration, respectively. Initial fracture tractions are determined according to the S-Lin-Lin fracture deformation model. The boundary conditions consist of combinations of roller, absorbing and time-dependent Dirichlet conditions and are illustrated in Figure~\ref{fig:bc_type_comp_shear}. 

We performed separate temporal and spatial convergence analyses with a wave amplitude of $U=50.0\,\mu\mathrm{m}$ and a frequency of $100.0\,\mathrm{kHz}$. Four refinement levels, $r=0,1,2,3$, were considered and used to determine the grid and time step size. An overview of the grid and time step size is shown in Table~\ref{tab:convergence_setup_4_1}. Note that the spatial analysis uses a finer time step than the finest refinement in the temporal analysis, which is done to ensure that the temporal error did not influence the spatial convergence.

\begin{table}[width=.9\linewidth,cols=4,pos=htbp]
    \caption{Grid and time step sizes used for the temporal and spatial convergence analyses for the linear spring-type deformation model}
    \label{tab:convergence_setup_4_1}
    \begin{tabular*}{\tblwidth}{@{} LLLL@{} }
    \toprule
    Analysis & Refinement levels & $\Delta x$ [mm] & $\Delta t$ [$\mu\mathrm{s}$]\\
    \midrule
    Temporal & $r=0,1,2,3$ & 0.125 & $1.0/2^r$\\
    Spatial & $r=0,1,2,3$ & $1.0/2^r$ & 0.03125\\
    \bottomrule
    \end{tabular*}
\end{table}

Convergence is assessed by examining the numerical errors in discrete solutions for bulk rock and fracture quantities. In the bulk rock, the discrete displacement $\boldsymbol u$ and traction $\boldsymbol q$ are considered, while on the fracture, the discrete displacement jump $\llbracket \boldsymbol u \rrbracket$ and fracture contact traction $\boldsymbol q_l$ are considered. Here and throughout the paper, bold symbols denote discrete quantities, distinguishing them from their continuous counterparts. Errors are computed using Equation~\ref{eq:error_cell} for cell quantities ($\boldsymbol u$, $\llbracket \boldsymbol u \rrbracket$, $\boldsymbol q_l$) and Equation~\ref{eq:error_face} for cell-face quantities ($\boldsymbol q$), all evaluated at the final time $t_{\mathrm{final}}=20.0\,\mu\mathrm{s}$. The relative errors are plotted against a spatial or temporal quantity which doubles for each refinement level. In the rock domain, spatial refinement quadruples the number of cells $N_x$, so $(N_x)^{1/2}$ is used as the $x$-axis quantity. Along the fracture, the number of cells $N_x$ doubles with each refinement, so $N_x$ is used as the $x$-axis quantity. For the temporal analysis, the number of time steps $N_t$ is used.

\begin{figure}[pos=htbp]
    \centering
    \includegraphics[width=0.95\textwidth]{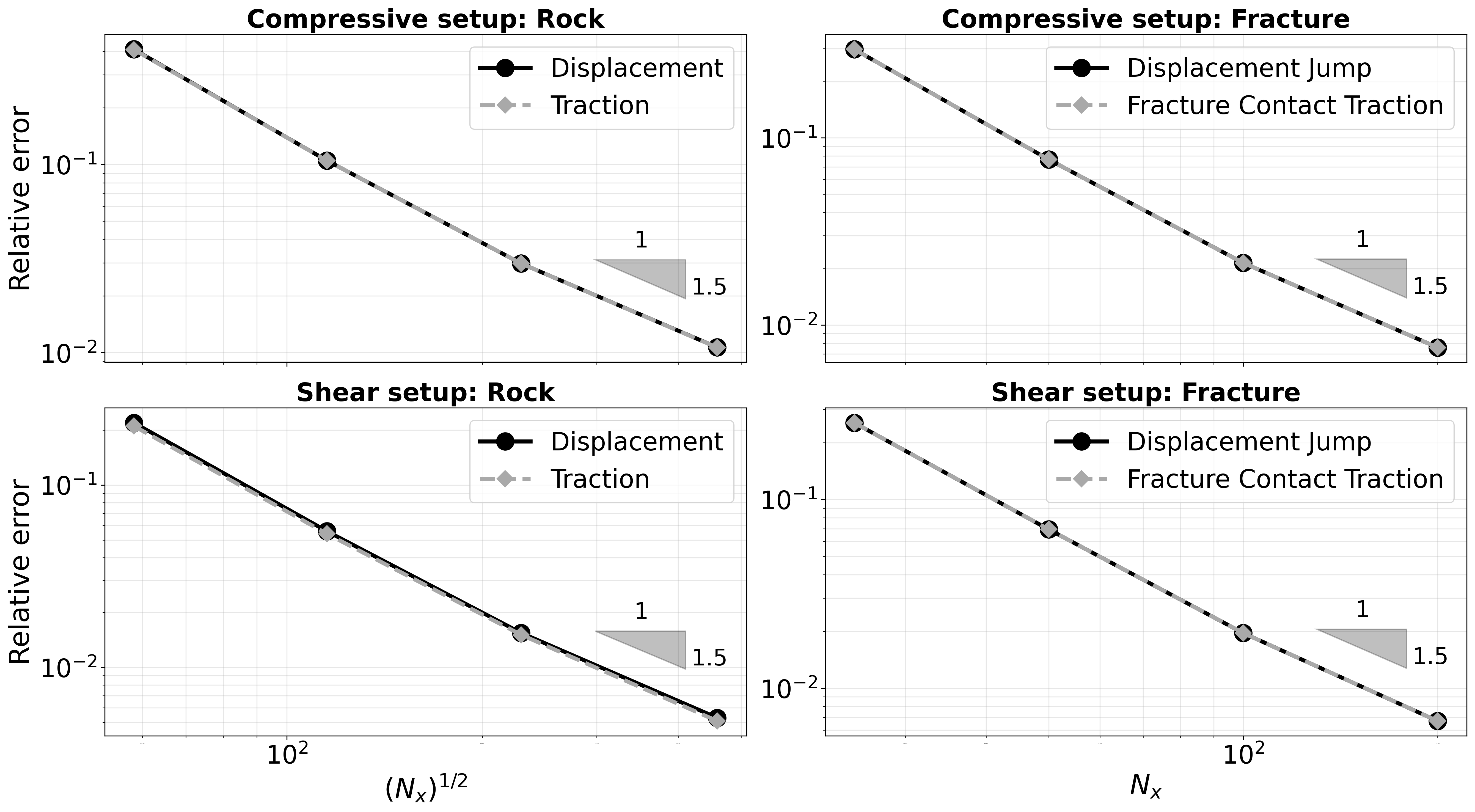}
    \caption{Spatial convergence analysis results for the S-Lin-Lin fracture deformation model. The relative errors in the rock displacement $\boldsymbol{u}$, cell-face traction $\boldsymbol{q}$, fracture displacement jump $\llbracket \boldsymbol{u} \rrbracket$ and fracture contact traction $\boldsymbol{q}_l$ converge with an order of 1.5. Errors are reported for both compressive and shear waves and are computed using Equations~\ref{eq:error_cell} and \ref{eq:error_face}}
    \label{fig:convergence_linear_spring_space}
\end{figure}

\begin{figure}[pos=htbp]
    \centering
    \includegraphics[width=0.95\textwidth]{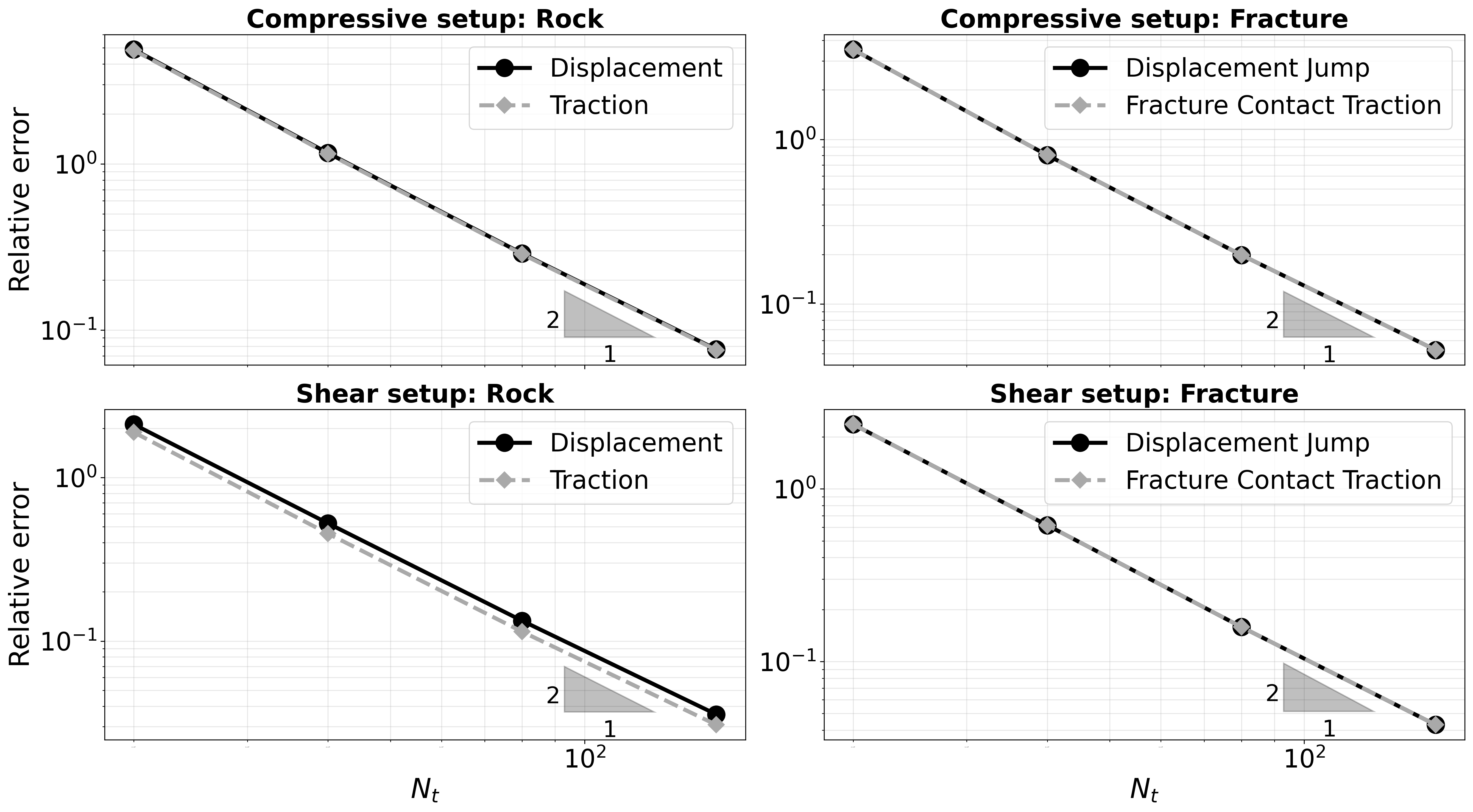}
    \caption{Temporal convergence analysis results for the S-Lin-Lin fracture deformation model. The relative errors in the rock displacement $\boldsymbol{u}$, cell-face traction $\boldsymbol{q}$, fracture displacement jump $\llbracket \boldsymbol{u} \rrbracket$ and fracture contact traction $\boldsymbol{q}_l$ converge with second-order accuracy. Errors are reported for both compressive and shear waves and are computed using Equations~\ref{eq:error_cell} and \ref{eq:error_face}}
    \label{fig:convergence_linear_spring_time}
\end{figure}

The results shown in Figures~\ref{fig:convergence_linear_spring_space}~and~\ref{fig:convergence_linear_spring_time} demonstrate second-order convergence in time and approximately 1.5th-order convergence in space. The temporal convergence rate is optimal for the Newmark method, while the spatial convergence rate is within the expected range for the MPSA-W discretization \citep{Keilegavlen2015,Berge2020}.

\FloatBarrier
\subsection{Spring-Type Fracture Deformation: Barton--Bandis Model}
\label{subsec:spring_type_BB}
In this section, we investigate the transmission of P-waves through a fracture governed by the S-Lin-BB model. 

Using results from the literature, we obtain an approximation of the theoretical transmission coefficient $T_{th}$ for P-waves propagating through S-Lin-BB fractures. The quantity $T_{th}$ is then used to compute the relative error in the numerical transmission coefficient $T_{num}$:
\begin{equation}
    \mathcal{E}_{T}=\frac{\lvert T_{th} - T_{num}\rvert}{\lvert T_{th} \rvert}, \label{eq:T_error}
\end{equation}
where the expressions for $T_{th}$ and $T_{num}$ are presented below.

\subsubsection{Theoretical transmission coefficient}\label{subsubsec:theroetical_transmission_coefficient}
The fracture deformation model permits discontinuous displacement across the fracture and, consequently, discontinuous velocity. The fracture has two sides, where side $2$ is defined as the side furthest from the wave source. The theoretical velocity on fracture side $2$ is used to compute the theoretical transmission coefficient $T_{th}$.

In \citet{Zhao2001}, a discrete evolution equation is derived for the theoretical particle velocity on side $2$ of a fracture governed by the S-Lin-BB model. The equation reads
\begin{equation}
    \dot u_{th, 2}(t_{i+1}) = 
    \frac{
    T_{inc} \left( 2\dot u_{th, 2}(t_i) - 2p\left(t_i - \frac{l}{c_p}\right) \right)
    }{
    M \left(
    -\frac{\rho c_p}{K_n + \frac{\rho c_p \dot u_{th, 2}(t_i)}{\Delta u_{\max}}}
    +
    \frac{
    (\rho c_p)^2 \dot u_{th, 2}(t_i)
    }{
    \Delta u_{\max}
    \left(
    K_n + \frac{\rho c_p \dot u_{th, 2}(t_i)}{\Delta u_{\max}}
    \right)^2
    }
    \right)
    }
    + \dot u_{th, 2}(t_i),
    \label{eq:velocity_barton_bandis}
\end{equation}
where $p(t)$ for $t \in [0,T_{inc}]$ denotes the incident particle velocity driving the transmission problem, $T_{inc}$ is the duration of the incident pulse and $M$ denotes the number of time steps within one period of the incident wave. Finally, $t_i$ for $i=0,\ldots,N_t$ denotes the discrete time levels spanning the interval $[0,t_{final}]$. Given the initial velocity on fracture side $2$, Equation~\ref{eq:velocity_barton_bandis} can be applied iteratively to compute $\dot u_{th,2}(t_{i+1})$ from $\dot u_{th,2}(t_i)$ up to $t_{final}$. The theoretical transmission coefficient is then computed as
\begin{equation} 
    T_{th}=\frac{\max(\dot u_{th,2}(l,t))}{\max(p(t))}.
\end{equation}

In this study, we have used $M=10^5$ time steps per period of the incident wave to ensure an approximation error of $T_{th}$ wihch is well below the discretization error.

\subsubsection{Numerical transmission coefficient and convergence results}\label{subsubsec:numerical_transmission_coefficient_and_convergence}
We consider a rectangular domain with a vertical fracture located at $x=l=25.0\,\mathrm{mm}$. The domain is similar to that shown in Figure~\ref{fig:domain_A_B}, with media A and B assigned identical material parameters, resulting in a homogeneous medium.
 
The system is initially at rest, and the wave is imposed through the boundary conditions. The boundary conditions are shown in Figure~\ref{fig:bc_type_comp_shear} with $u_x = U\sin^4(\omega t)$ and $u_y=0$ for $t \in [0,5.0\cdot 10^{-6}]$. This boundary condition corresponds to the incident velocity wave $\dot u(0, t)=\left[p(t),\;0\right]^T$ with $p(t)=4U\omega\sin^3(\omega t)\cos(\omega t)$. The convergence analysis is performed with displacement amplitude $U=5.0\cdot10^{-8}\,\mathrm{m}$, frequency $f=100.0\,\mathrm{kHz}$ and $\Delta u_{\max}=1.0\cdot10^{-7}\,\mathrm{m}$.

\begin{figure}[pos=htbp]
    \centering
    \includegraphics[width=0.6\textwidth]{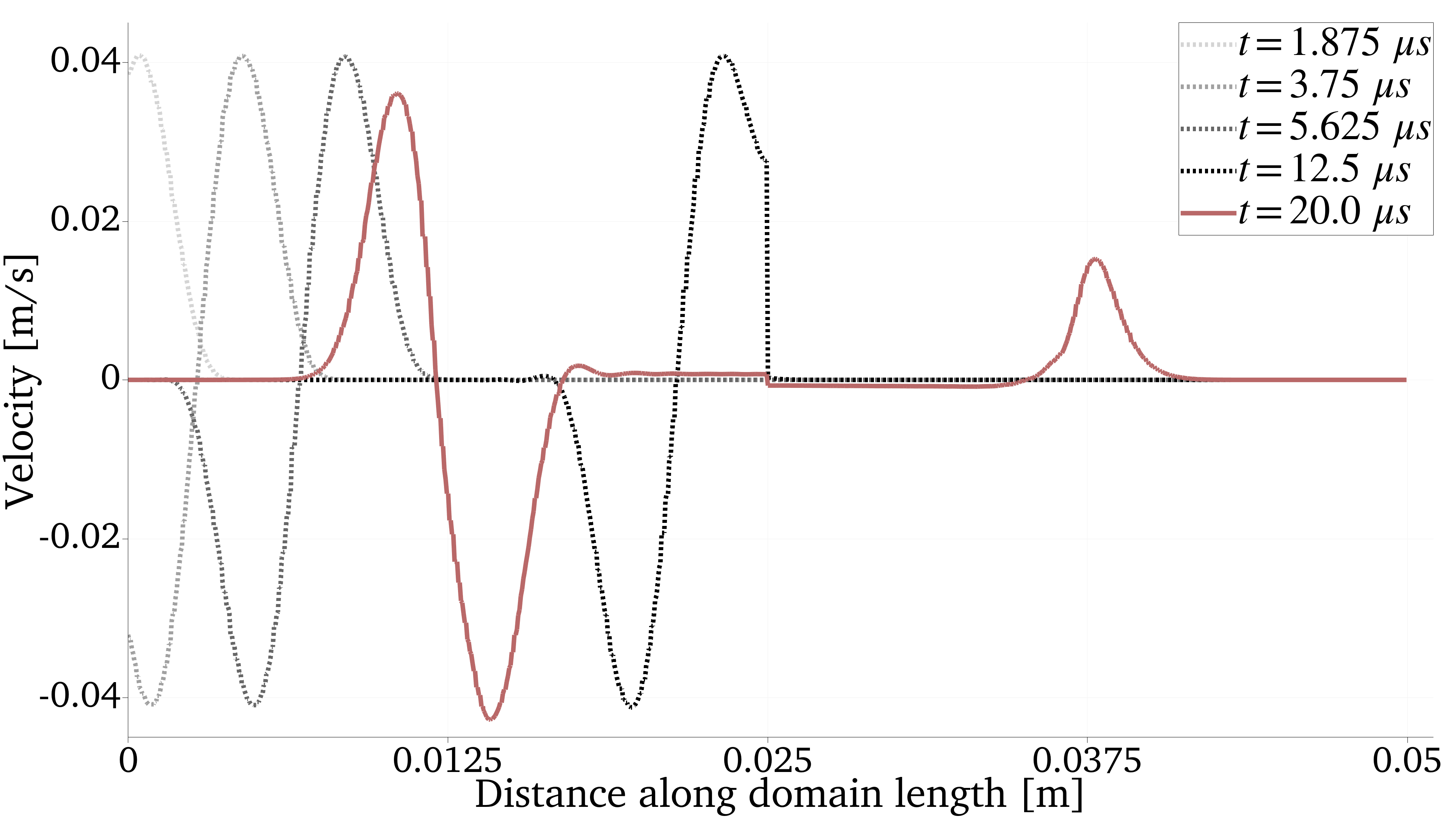}
    \caption{Velocity along the domain length at five different times from the convergence study using the S-Lin-BB model. The fracture is located at $x = 0.025\,\mathrm{m}$, and the gray and black (dotted) curves show the wave field before transmission through the fracture. The red (solid) curve represents the wave field at the final simulation time. The transmitted wave is located near $x = 0.0375\,\mathrm{m}$. The results are obtained from the finest refinement, namely $\Delta x=0.125\,\mathrm{mm}$ and $\Delta t = 0.0625\,\mu \mathrm{s}$}
    \label{fig:spring_BB_convergence}
\end{figure}

Four refinement levels, $r=0,1,2,3$, were considered and used to determine the grid size and time step size. An overview of the grid and time step size is shown in Table~\ref{tab:convergence_setup_4_2}. The simulation is performed up to the final time $t_{final}=20.0\,\mu\mathrm{s}$, which allows sufficient time for the wave peak to propagate through the fracture. With reference to Figure~\ref{fig:domain_A_B}, the numerical transmission coefficient is defined as
\begin{equation}
    T_{num} = \frac{\max(\dot{\boldsymbol{u}}_B)}{\max(p(t))},
\end{equation}
where $\dot{\boldsymbol{u}}_B$ denotes the particle velocity in region B. Figure~\ref{fig:spring_BB_convergence} illustrates the velocity field at five different times during the simulation. The quantity $\max(\dot{\boldsymbol{u}}_B)$ corresponds to the peak of the red solid-line wave located to the right of the fracture.

\begin{table}[width=.9\linewidth,cols=4,pos=htbp]
    \caption{Grid and time step sizes used for the temporal and spatial convergence analyses for the Barton--Bandis spring-type deformation model}
    \label{tab:convergence_setup_4_2}
    \begin{tabular*}{\tblwidth}{@{} LLLL@{} }
    \toprule
    Analysis & Refinement levels & $\Delta x$ [mm] & $\Delta t$ [$\mu\mathrm{s}$]\\
    \midrule
    Temporal & $r=0,1,2,3$ & $0.125$ & $0.5/2^r$\\
    Spatial & $r=0,1,2,3$ & $1.0/2^r$ & $0.0625$\\
    \bottomrule
    \end{tabular*}
\end{table}

\begin{figure}[pos=htbp]
    \centering
    \includegraphics[width=0.5\textwidth]{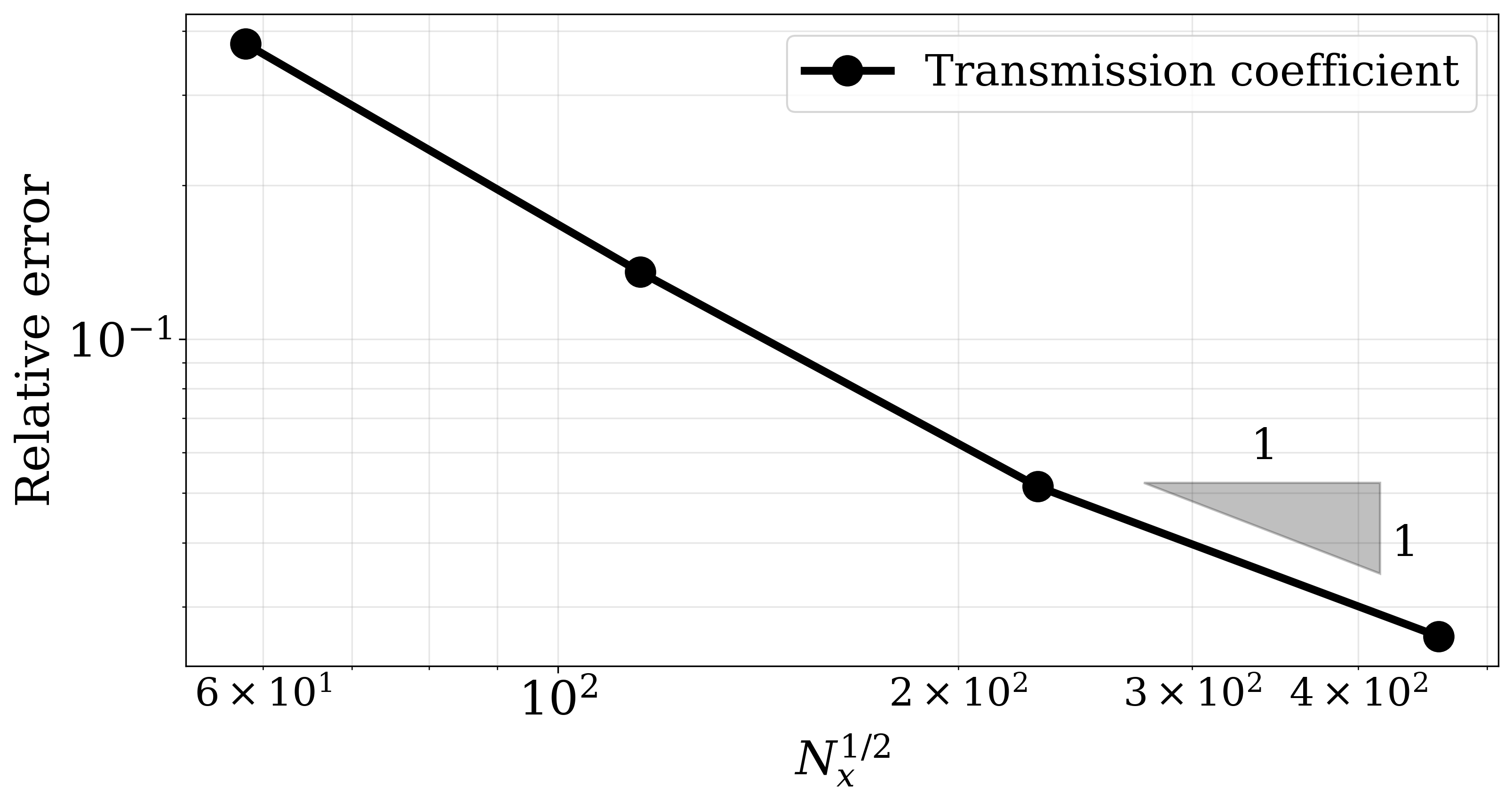}
    \caption{Spatial convergence analysis for the S-Lin-BB fracture deformation model. The plot shows that relative error in the numerical transmission coefficient $T_{num}$ converges with order 1. The error is computed by Equation~\ref{eq:T_error}} 
    \label{fig:convergence_nonlinear_spring_space}
\end{figure}

\begin{figure}[pos=htbp]
    \centering
    \includegraphics[width=0.5\textwidth]{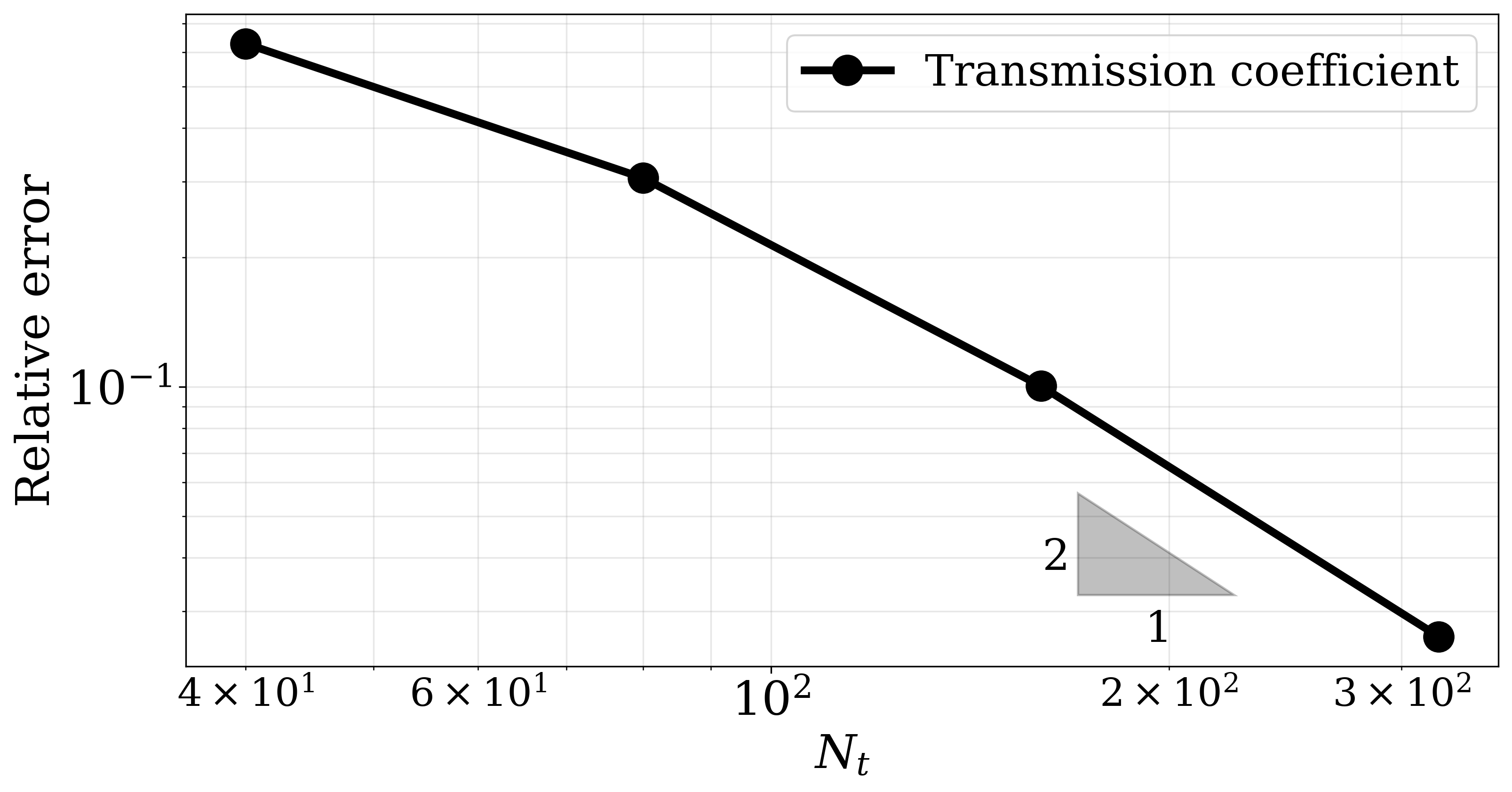}
    \caption{Temporal convergence analysis for the S-Lin-BB fracture deformation model. The plot shows that relative error in the numerical transmission coefficient $T_{num}$ converges with order 2. The error is computed by Equation~\ref{eq:T_error}} 
    \label{fig:convergence_nonlinear_spring_time}
\end{figure}

The convergence results for the numerical transmission coefficient, shown in Figures~\ref{fig:convergence_nonlinear_spring_space}~and~\ref{fig:convergence_nonlinear_spring_time}, indicate first order convergence in space and second order convergence in time. The spatial convergence is lower than that observed in the previous section for displacements and tractions, which were measured in an $L^2$-type norm. Since the transmission operator $T_{\mathrm{num}}$ has a more $L^\infty$-like character, convergence rates exceeding second order cannot generally be expected in this setting, and second-order accuracy is therefore considered optimal. These results are consistent with numerical findings for multi-point flux approximation methods, which are scalar analogues of the MPSA-W method considered herein. In particular, convergence rates in $L^\infty$-type norms were found to be only marginally lower than those obtained in $L^2$-type norms for certain test cases \citep{Aavatsmark2008}. Similar results for the MPSA-W method have not previously been reported.

\FloatBarrier
\subsection{Fracture Contact Mechanics with Friction}
\label{subsec:fracture_contact_mechanics_friction}
In this section, we perform self-convergence studies of the C-Coul-Lin and C-Coul-BB contact models, extending the previous spring-based analyses to include contact and friction. For both models, convergence is assessed for the rock displacement, fracture contact traction and displacement jump across the fracture following the procedure described in Appendix~\ref{appendix_nonmatching_grids}.

\begin{figure}[pos=htbp]
    \centering
    \includegraphics[width=0.6\textwidth]{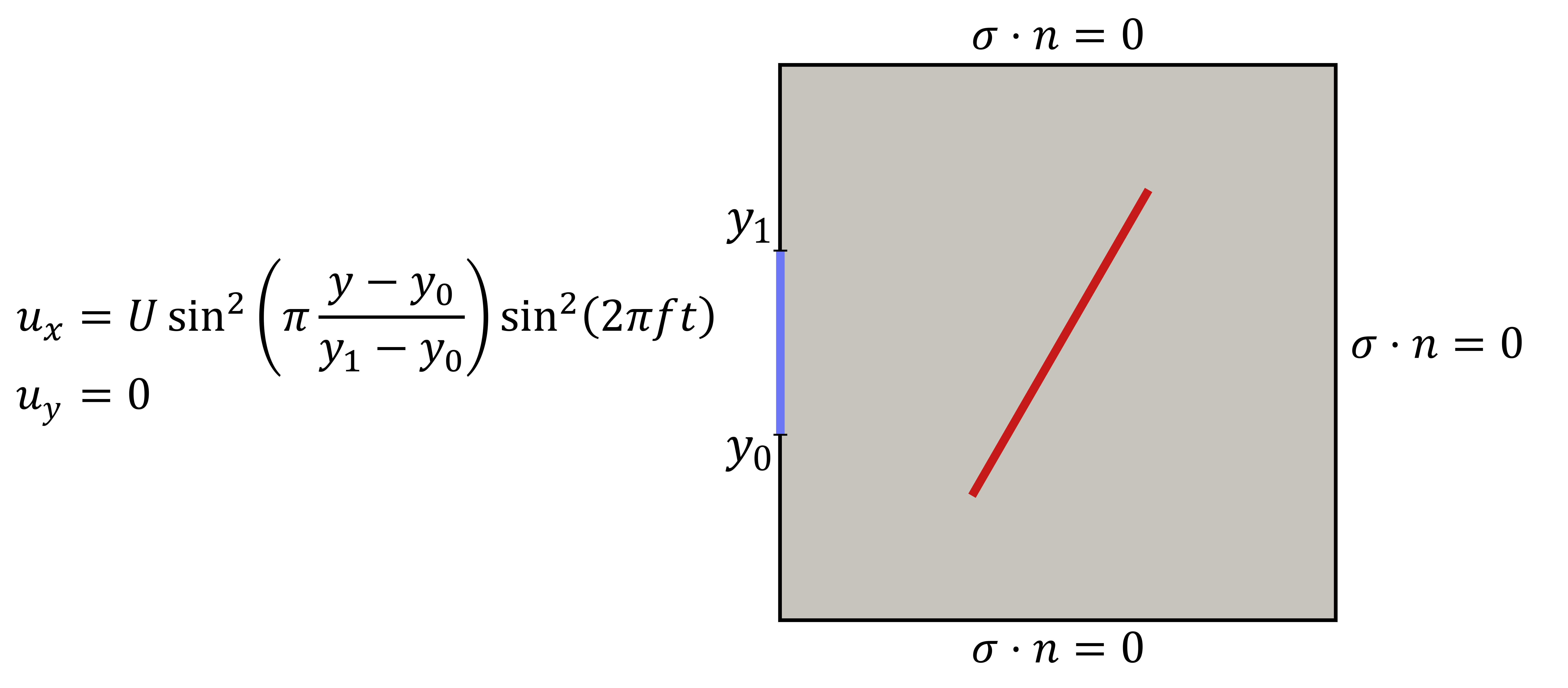}
    \caption{Simulation domain and boundary conditions for the self-convergence analyses with the C-Coul-Lin and C-Coul-BB models. The Dirichlet condition on the left boundary imposes the wave and is nonzero only between $y_0$ and $y_1$. The remaining boundaries are stress-free}
    \label{fig:bc_type_self_convergence}
\end{figure}

The simulation domain is two-dimensional and measures $25.0\,\mathrm{mm}$ by $25.0\,\mathrm{mm}$. The domain contains a $16.0\,\mathrm{mm}$ long diagonal fracture with friction coefficient $F=1.0$ and maximum fracture closure $\Delta u_{\max} = 5.0\cdot 10^{-5}\,\mathrm{m}$. Stress-free conditions are prescribed on the right, top and bottom boundaries, while a tapered time-dependent Dirichlet condition is imposed in the $x$-direction on the left boundary:
\begin{align}
    u(0,y,t) =
    \begin{cases}
        \begin{bmatrix}
        U \sin^2\!\left(\pi \frac{y - y_0}{y_1 - y_0}\right)\sin^2(2\pi f t),\; 0
        \end{bmatrix}^{T}
        & \text{for } y_0 \le y \le y_1, \\
        \begin{bmatrix}
        0,\; 0
        \end{bmatrix}^{T}
        & \text{otherwise,}
    \end{cases}
    \label{eq:dirichlet_tapered_condition}
\end{align}
where $y_0=6.25\,\mathrm{mm}$ and $y_1=18.75\,\mathrm{mm}$. Figure~\ref{fig:bc_type_self_convergence} shows a schematic of the fractured domain and the associated boundary conditions.

The simulation is performed up to the final time $t_{final}=8.75\,\mu\mathrm{s}$. The refinement level $r$ determines the grid size and time step size, and an overview of the different values are shown in Table~\ref{tab:convergence_setup_4_3}. We use the solution with the highest $r$ as the reference solution.

\begin{table}[width=.9\linewidth,cols=4,pos=htbp]
    \caption{Grid and time step sizes used for the temporal and spatial convergence analyses for the contact mechanics fracture deformation models}
    \label{tab:convergence_setup_4_3}
    \begin{tabular*}{\tblwidth}{@{} LLLL@{} }
    \toprule
    Analysis & Refinement levels & $\Delta x$ [mm] & $\Delta t$ [$\mu\mathrm{s}$]\\
    \midrule
    Temporal & $r=0,1,2,3,4,5$ & $0.1$ & $0.21875/2^r$\\
    Spatial & $r=0,1,2,3,4$ & $1.6/2^r$ & $0.21875/2^5$\\
    \bottomrule
    \end{tabular*}
\end{table}

\begin{figure}[pos=htbp]
    \centering
    \includegraphics[width=0.8\textwidth]{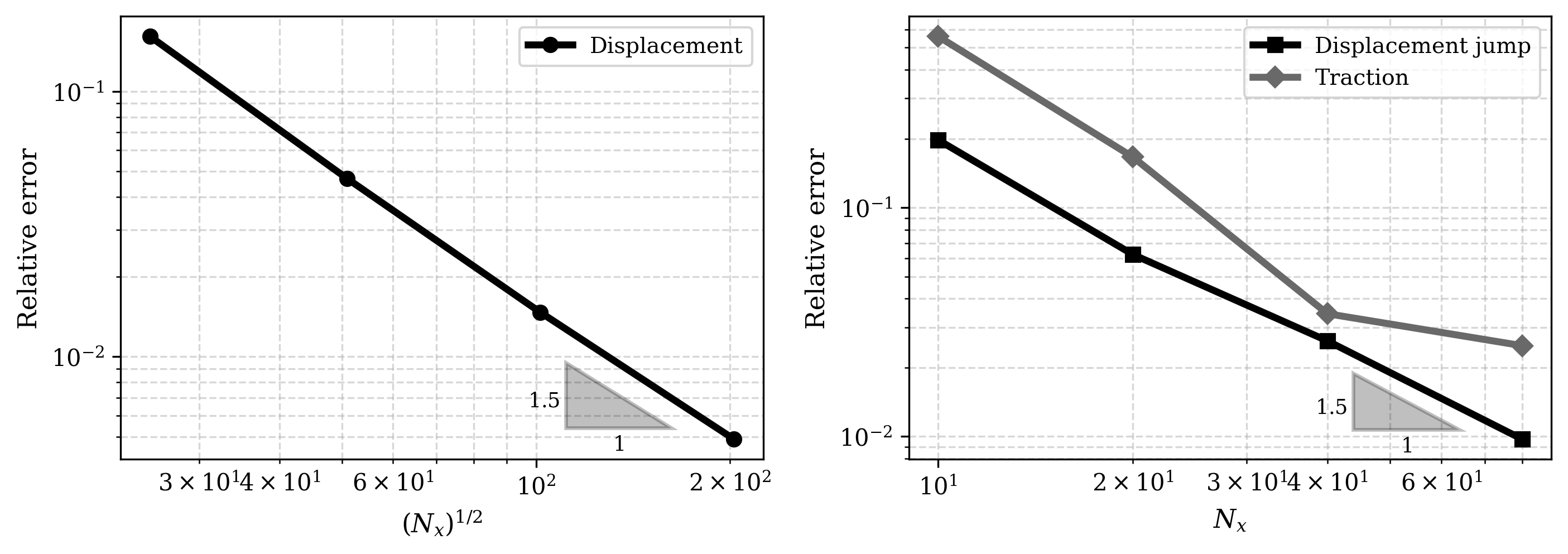}
    \caption{Spatial convergence analysis for the C-Coul-Lin fracture deformation model. Left: Relative error in the rock displacement $\boldsymbol{u}$. Right: Relative errors in fracture displacement jump $\llbracket \boldsymbol{u} \rrbracket$ and fracture contact traction $\boldsymbol{q}_l$. All quantities converge with an order of approximately 1.5. The errors are computed using the procedure described in Appendix~\ref{appendix_nonmatching_grids}}
    \label{fig:self_convergence_LM_space}
\end{figure}

\begin{figure}[pos=htbp]
    \centering
    \includegraphics[width=0.8\textwidth]{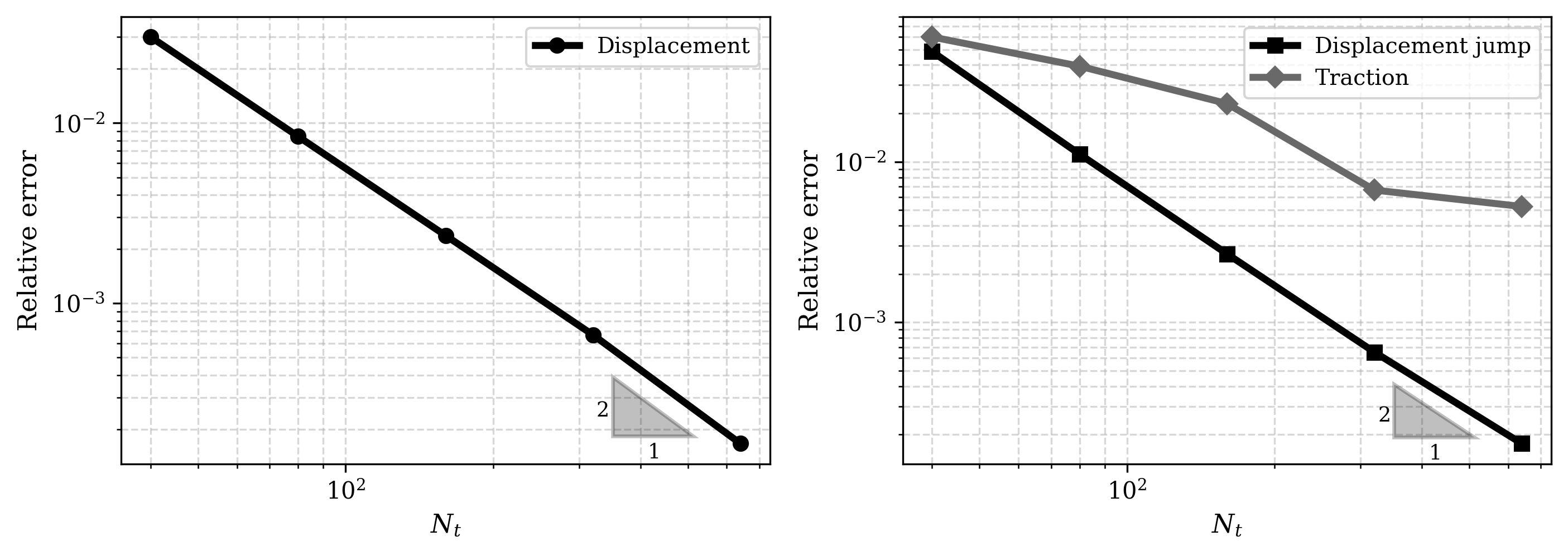}
    \caption{Temporal convergence analysis for the C-Coul-Lin fracture deformation model. Left: Relative error in the rock displacement $\boldsymbol{u}$. Right: Relative errors in fracture displacement jump $\llbracket \boldsymbol{u} \rrbracket$ and fracture contact traction $\boldsymbol{q}_l$. The displacement unknowns converge with an order of 2, while the fracture contract traction does not exceed a convergence order of 1.5. All refinements of the solution here live on the same grid, and the errors are thus computed as regular L2-type errors without the need for the interpolation and transfer grids described in  Appendix~\ref{appendix_nonmatching_grids}}
    \label{fig:self_convergence_LM_time}
\end{figure}

\begin{figure}[pos=htbp]
    \centering
    \includegraphics[width=0.8\textwidth]{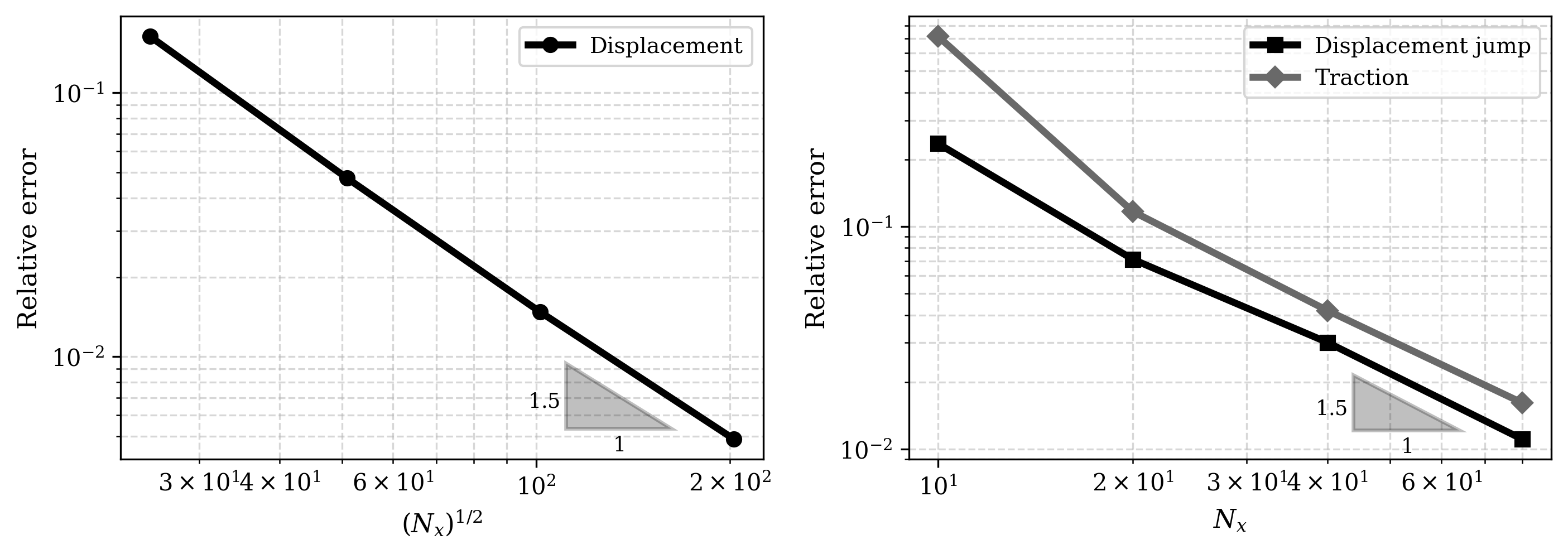}
    \caption{Spatial convergence analysis for the C-Coul-BB fracture deformation model. Left: Relative error in the rock displacement $\boldsymbol{u}$. Right: Relative errors in fracture displacement jump $\llbracket \boldsymbol{u} \rrbracket$ and fracture contact traction $\boldsymbol{q}_l$. All quantities converge with an order of approximately 1.5. The errors are computed using the procedure described in Appendix~\ref{appendix_nonmatching_grids}}
    \label{fig:self_convergence_BB_space}
\end{figure}

\begin{figure}[pos=htbp]
    \centering
    \includegraphics[width=0.8\textwidth]{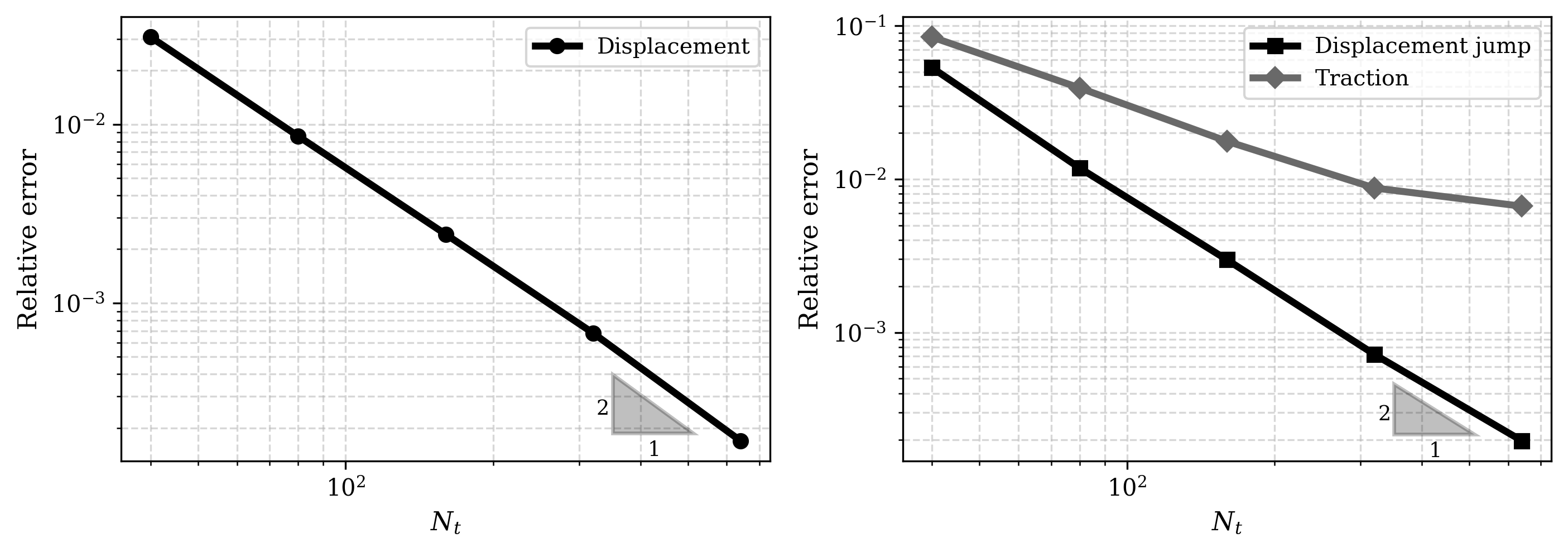}
    \caption{Temporal convergence analysis for the C-Coul-BB fracture deformation model. Left: Relative error in the rock displacement $\boldsymbol{u}$. Right: Relative errors in fracture displacement jump $\llbracket \boldsymbol{u} \rrbracket$ and fracture contact traction $\boldsymbol{q}_l$. The displacement unknowns converge with an order of 2, while the fracture contract traction does not exceed a convergence order of 1.5. All refinements of the solution here live on the same grid, and the errors are thus computed as regular L2-type errors without the need for the interpolation and transfer grids described in  Appendix~\ref{appendix_nonmatching_grids}}
    \label{fig:self_convergence_BB_time}
\end{figure}

Referring to Figures~\ref{fig:self_convergence_LM_space} and \ref{fig:self_convergence_BB_space}, the spatial errors in bulk displacement, displacement jump and fracture contact traction converge with orders between one and two. The temporal convergence results shown in Figures~\ref{fig:self_convergence_LM_time}~and~\ref{fig:self_convergence_BB_time} exhibit second-order convergence for the displacement quantities, consistent with the theoretical accuracy of the Newmark method \citep{Subbaraj1989}, whereas the fracture contact traction converges at a lower rate. The convergence of the fracture contact traction may be influenced by the oscillatory behavior known to occur near fracture tips \citep{Berge2020}, and possibly also transitions between different contact regimes.

\begin{figure}[pos=htbp]
    \centering
    \includegraphics[width=0.6\textwidth]{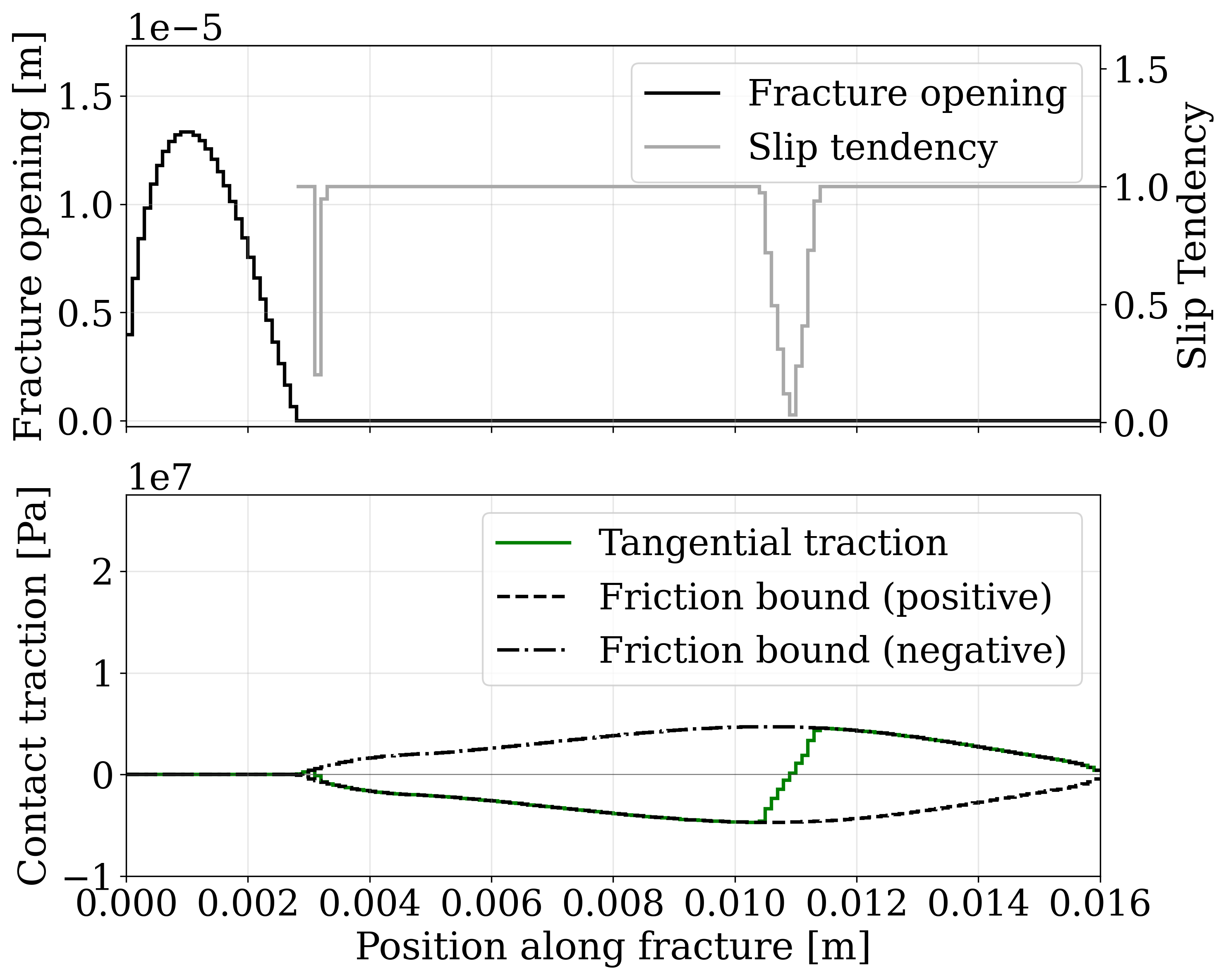}
    \caption{Fracture opening $\delta$, slip tendency $s$, tangential fracture traction $\boldsymbol{q}_\tau$ and friction bounds $\pm \boldsymbol{b}$ at the final time for the C-Coul-Lin model, using the reference solution for the spatial convergence analysis. Regions with $\delta > 0$ correspond to $s = \mathrm{NaN}$ and $\boldsymbol{q}_\tau = \boldsymbol{b} = 0$, regions with $s = 1$ correspond to $\boldsymbol{q}_\tau = \pm \boldsymbol{b}$ and regions with $s < 1$ correspond to $-\boldsymbol{b} < \boldsymbol{q}_\tau < \boldsymbol{b}$}
    \label{fig:lineplot_CL}
\end{figure}

\begin{figure}[pos=htbp]
    \centering
    \includegraphics[width=0.6\textwidth]{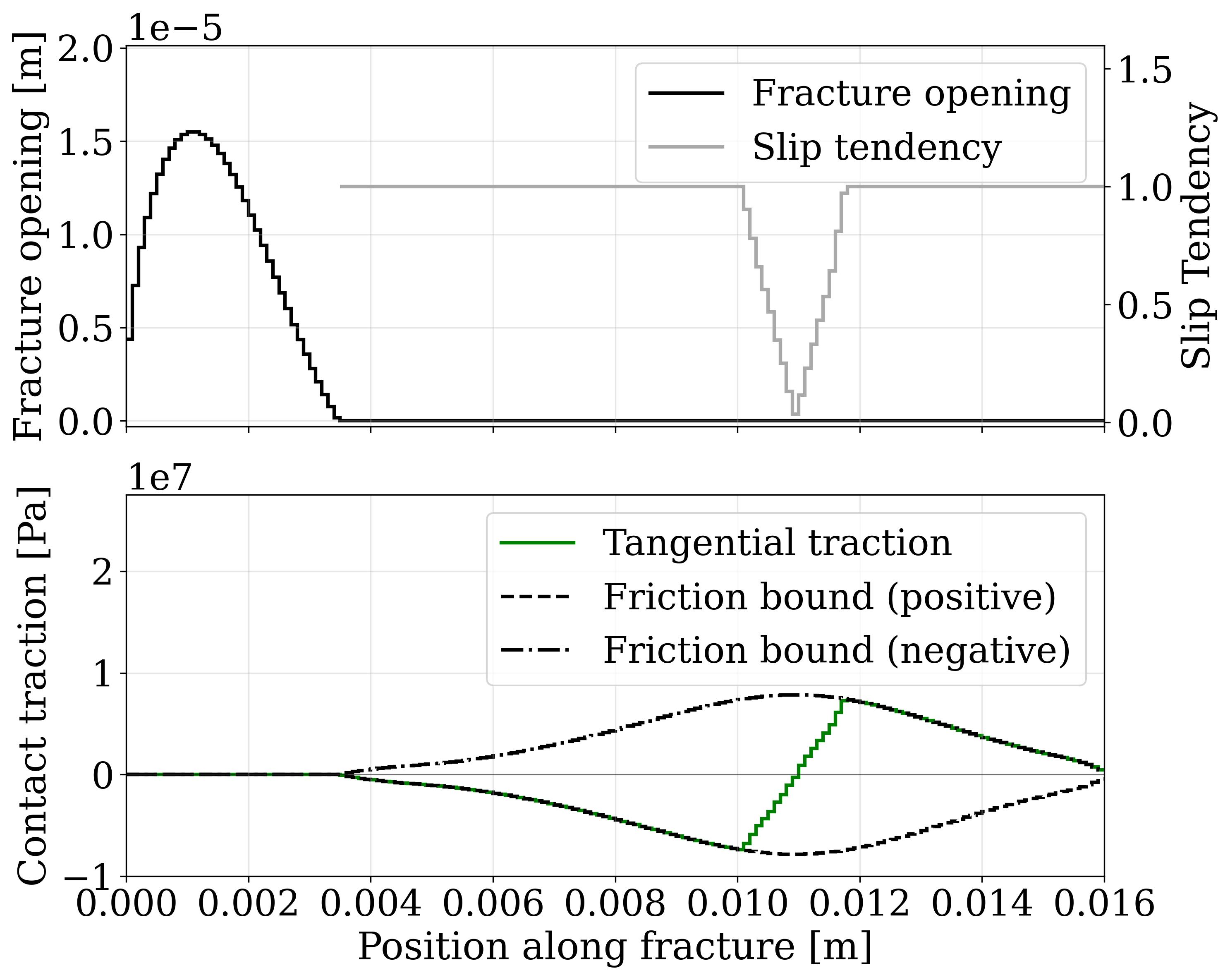}
    \caption{Fracture opening $\delta$, slip tendency $s$, tangential fracture traction $\boldsymbol{q}_\tau$ and friction bounds $\pm \boldsymbol{b}$ at the final time for the C-Coul-BB model, using the reference solution for the spatial convergence analysis. Regions with $\delta > 0$ correspond to $s = \mathrm{NaN}$ and $\boldsymbol{q}_\tau = \boldsymbol{b} = 0$, regions with $s = 1$ correspond to $\boldsymbol{q}_\tau = \pm \boldsymbol{b}$ and regions with $s < 1$ correspond to $-\boldsymbol{b} < \boldsymbol{q}_\tau < \boldsymbol{b}$}
    \label{fig:lineplot_CBB}
\end{figure}

\begin{figure}[pos=htbp]
    \centering
    \includegraphics[width=0.7\textwidth]{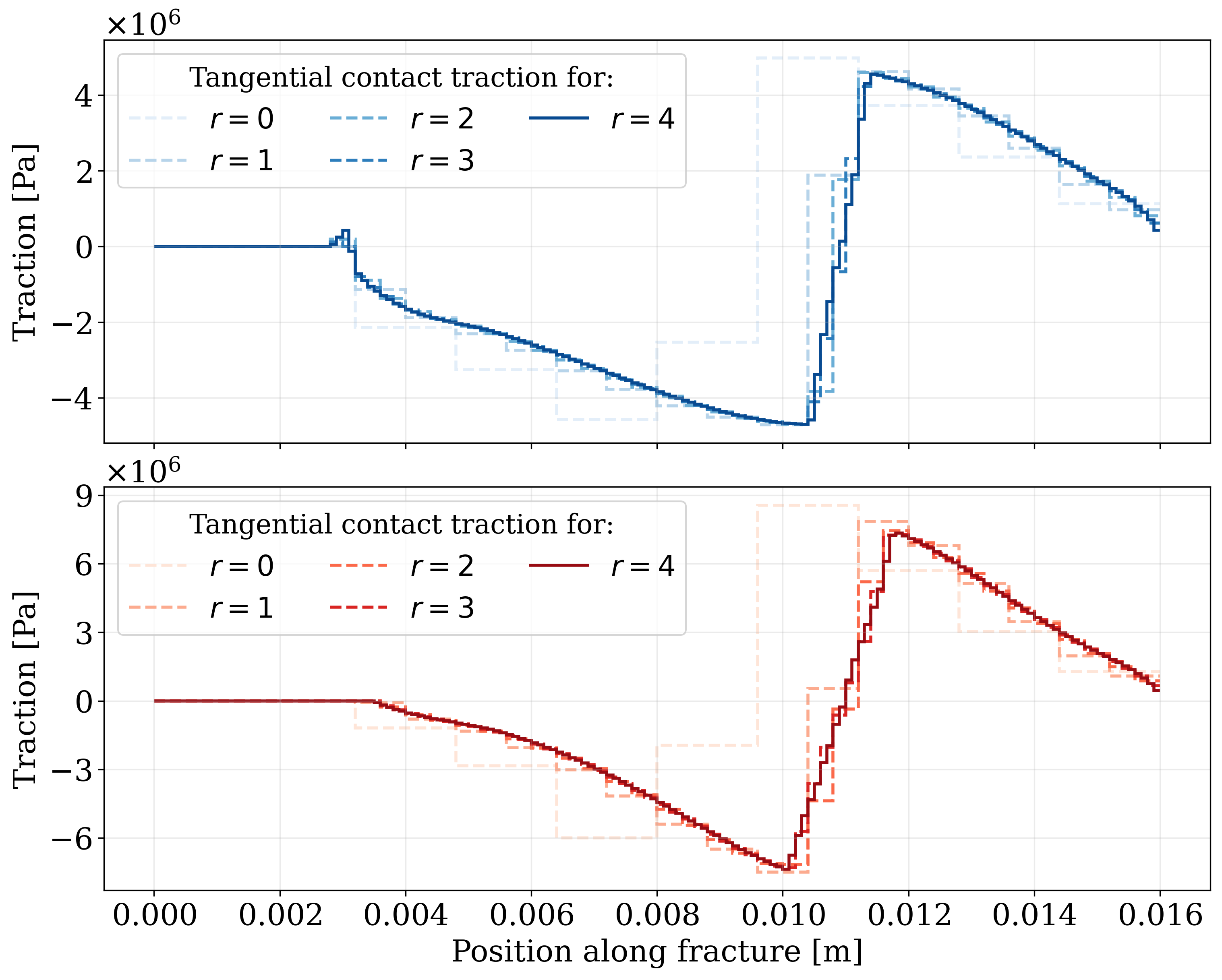}
    \caption{Tangential fracture contact traction for all refinements $r$ in the spatial convergence analysis, shown for C-Coul-Lin (top, blue) and C-Coul-BB (bottom, red). The coarser solutions are represented by dashed lines, while the reference solution is represented by a solid line}
    \label{fig:contact_traction_all_cell_sizes}
\end{figure}

\begin{figure}[pos=htbp]
    \centering
    \includegraphics[width=0.7\textwidth]{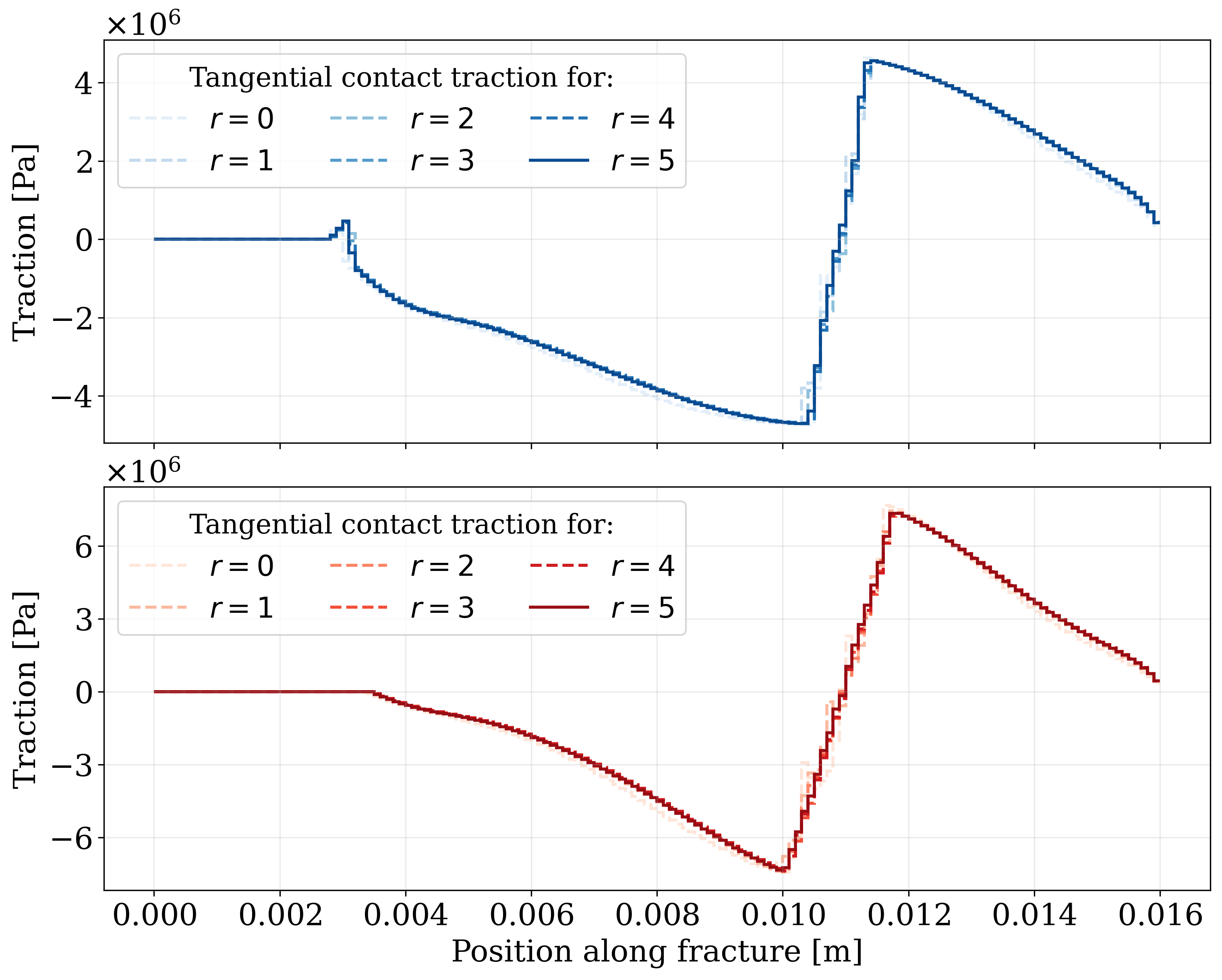}
    \caption{Tangential fracture contact traction for all refinements $r$ in the temporal convergence analysis, shown for C-Coul-Lin (top, blue) and C-Coul-BB (bottom, red). The coarser solutions are represented by dashed lines, while the reference solution is represented by a solid line}
    \label{fig:contact_traction_all_time_step_sizes}
\end{figure}

The convergence results are evaluated at a time when all three contact states are present along the fracture: open, sticking and sliding. This is illustrated in Figures~\ref{fig:lineplot_CL} and~\ref{fig:lineplot_CBB} for the C-Coul-Lin and C-Coul-BB models, respectively. Figures~\ref{fig:contact_traction_all_cell_sizes} and~\ref{fig:contact_traction_all_time_step_sizes} show the tangential fracture contact traction for all grid and time step sizes used in the self-convergence analyses. As the discretization is refined, the solutions approach the reference solution, consistent with the convergence behavior observed in the quantitative error measures. Overall, the results indicate stable and convergent behavior across all three contact regimes.

\FloatBarrier
\section{Comparison of the Fracture Deformation Models}\label{sec:comparison_of_fracture_deformation_models}
In this section, we compare simulation results from four different fracture deformation models S-Lin-Lin, S-Lin-BB, C-Coul-Lin and C-Coul-BB. For all four models, we present the resulting wave fields in the bulk rock, as well as the normal and tangential fracture deformation. Experimental studies have shown that fracture deformation is nonlinear under sufficiently large excitation  \citep{Bandis1983,Zhao2008Dynamic}, and the hyperbolic Barton--Bandis relation is a commonly used description of such nonlinear fracture deformation \citep{Bandis1983}. We therefore include both linear and nonlinear fracture deformation models in the comparison, allowing us to examine how the choice of deformation model affects the fracture deformation and the resulting wave field. The fracture parameters are chosen as representative values for the problem and are not calibrated against laboratory measurements. The results should therefore be interpreted as an illustration of the different bulk and fracture deformation responses captured by the four models, rather than as a quantitative assessment or validation of the models.

The numerical setup is identical to that presented in Section~\ref{subsec:fracture_contact_mechanics_friction}, except that $\Delta u_{\max} = 1.0\cdot10^{-5}\,\mathrm{m}$ and the final simulation time is increased to $t_{final}=1.05\cdot10^{-5}\,\mathrm{s}$. The longer simulation time ensures that the wave has fully propagated past the fracture before the results are analyzed. The simulation is performed with a time step of $\Delta t=5.25\cdot10^{-8}\,\mathrm{s}$ and a grid size of $\Delta x=0.2\,\mathrm{mm}$. Both the time step and grid size provide at least 10 points per wave period and wavelength, respectively, ensuring that the wave is adequately resolved. Also, the chosen grid and time step sizes are supported by the convergence analyses presented in Section~\ref{sec:numerical_convergence_analyses}, which show that they lie well within the convergent regimes for all the deformation models.

\begin{figure}[pos=ht]
    \centering
    \includegraphics[width=0.99\textwidth]{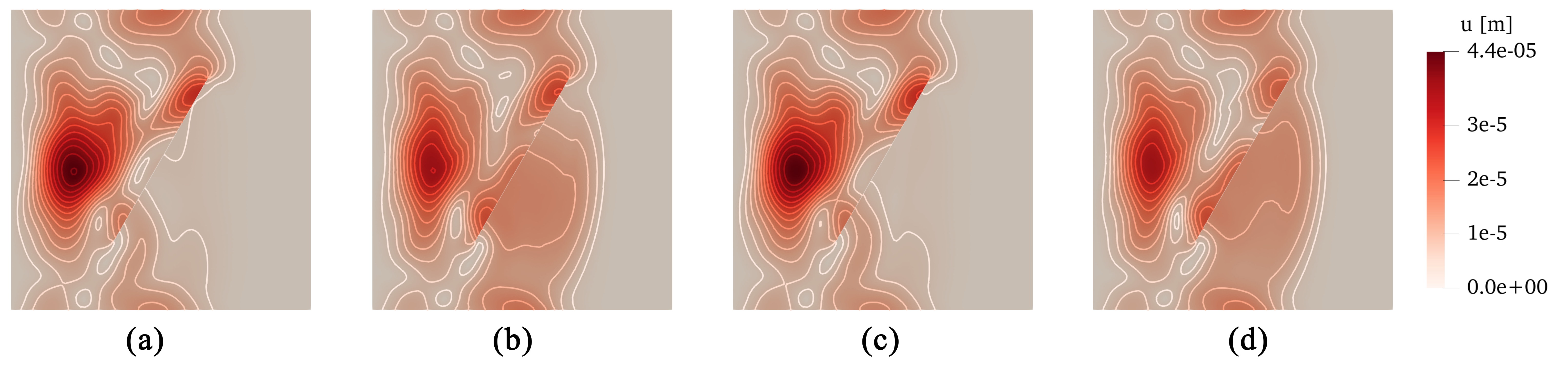}
    \caption{Magnitude of the displacement wave field for all models at the final time: (a) S-Lin-Lin, (b) S-Lin-BB, (c) C-Coul-Lin and (d) C-Coul-BB. Contour lines are added for visualization purposes and represent the same contour values in all four panels}
    \label{fig:wavefield_all_models}
\end{figure}
Figure~\ref{fig:wavefield_all_models} shows the wave fields obtained with the four fracture deformation models. The largest differences are observed in the transmitted wave field to the right of the fracture. Models with Barton--Bandis normal deformation (S-Lin-BB and C-Coul-BB; Figures~\ref{fig:wavefield_all_models} b) and d)) transmit significantly larger wave amplitudes than models with linear elastic normal deformation (S-Lin-Lin and C-Coul-Lin; Figures~\ref{fig:wavefield_all_models} a) and c)). While all models share the same normal stiffness $K_n$, the effective stiffness of the Barton--Bandis models increases during compression due to the $\Delta u_{\max}$ parameter, resulting in more efficient wave transmission. The wave fields predicted by the S-Lin-Lin and C-Coul-Lin models are very similar, differing only slightly in the wavefront contours, whereas those obtained with the S-Lin-BB and C-Coul-BB models exhibit noticeably different transmitted wavefront shapes. The primary difference between the S-Lin-BB and C-Coul-BB models is the width of the transmitted wave, which is greater for the contact model. Overall, the results indicate that the additional complexity introduced by contact mechanics and Barton--Bandis normal deformation can have a significant impact on the predicted wave field for the setup considered here.

\begin{figure}[pos=htbp]
    \centering
    \includegraphics[width=0.99\textwidth]{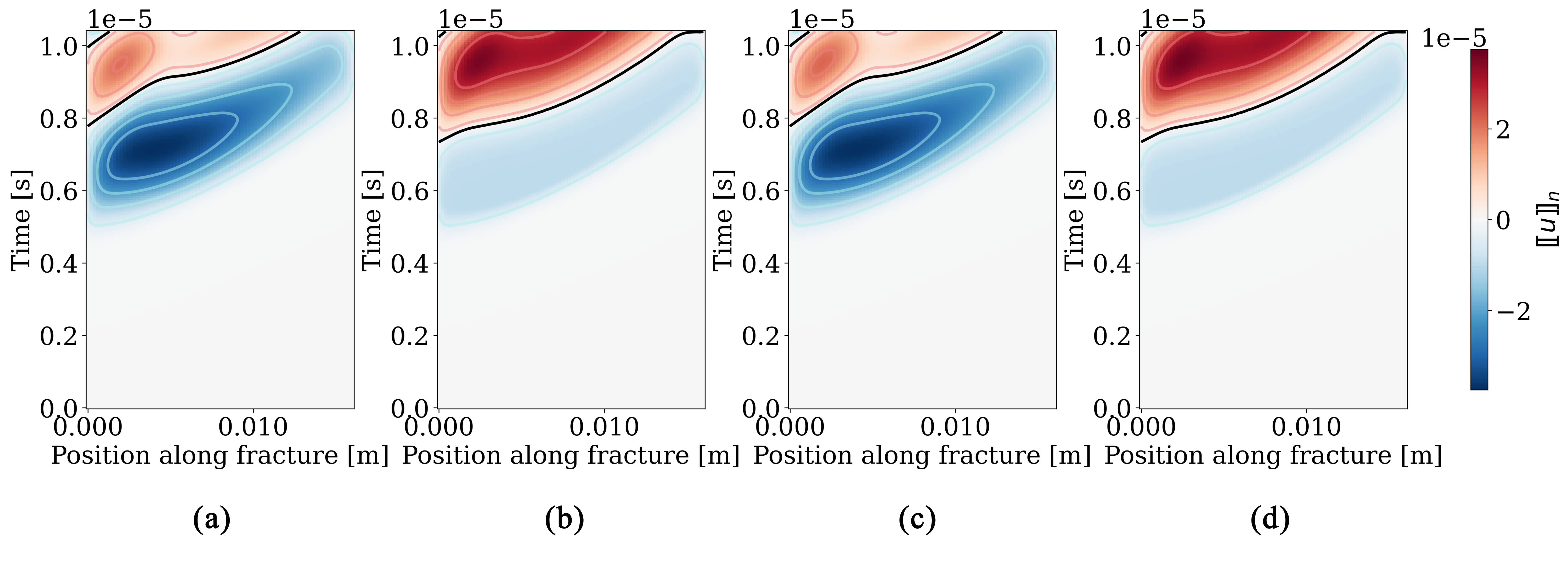}
    \caption{Normal fracture deformation $\llbracket \boldsymbol u \rrbracket_n$ for the entire simulation. (a) S-Lin-Lin, (b) S-Lin-BB, (c) C-Coul-Lin and (d) C-Coul-BB. Contour lines are added for visualization purposes and represent the same contour values in all four panels. The black contour line represents the zero contour}
    \label{fig:normal_all_models}
\end{figure}

\begin{figure}[pos=htbp]
    \centering
    \includegraphics[width=0.99\textwidth]{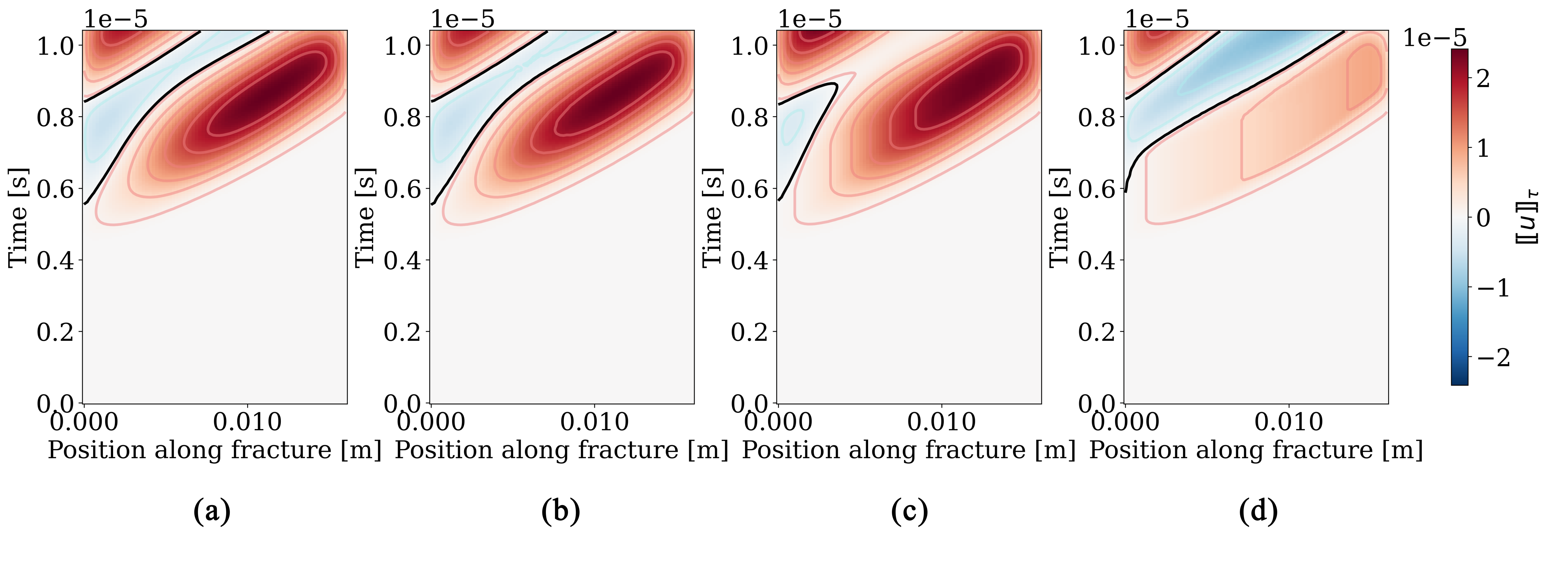}
    \caption{Tangential fracture deformation $\llbracket \boldsymbol u \rrbracket_\tau$ for the entire simulation. (a) S-Lin-Lin, (b) S-Lin-BB, (c) C-Coul-Lin and (d) C-Coul-BB. Contour lines are added for visualization purposes and represent the same contour values in all four panels. The black contour line represents the zero contour}
    \label{fig:tangential_all_models}
\end{figure}

\begin{figure}[pos=htbp]
    \centering
    \includegraphics[width=0.99\linewidth]{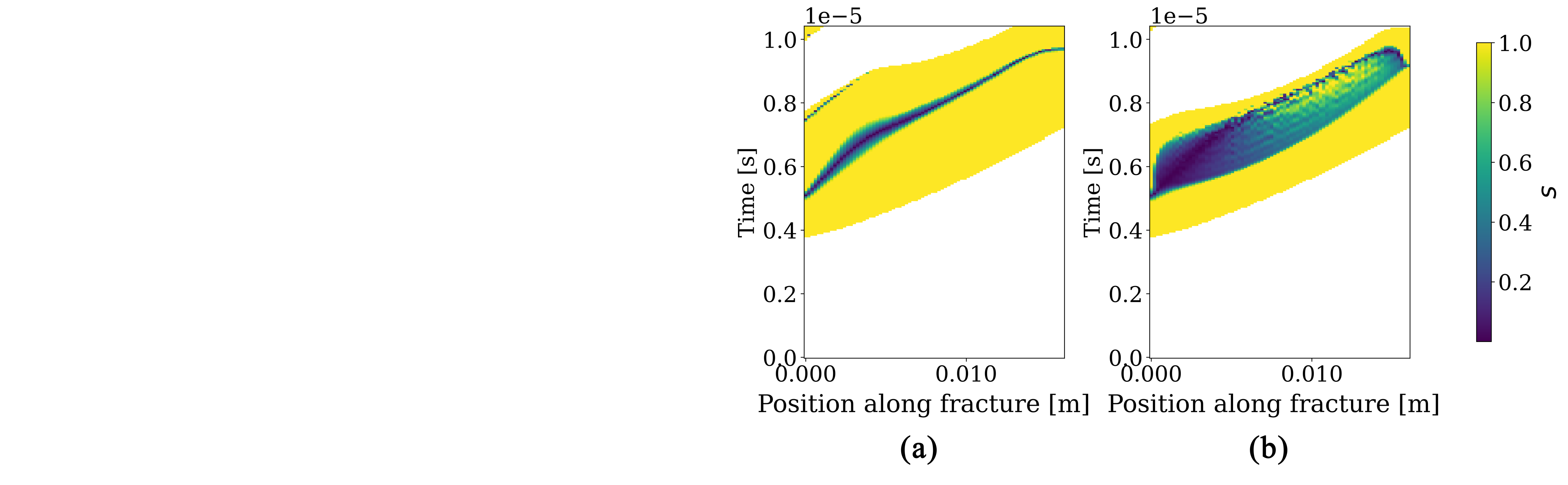}
    \caption{Slip tendency $s$ along the fractures during the entire simulation for the (a) C-Coul-Lin and (b) C-Coul-BB models.}
    \label{fig:slip_tendency_model_comparison}
\end{figure}

Figures~\ref{fig:normal_all_models} and~\ref{fig:tangential_all_models} show the normal and tangential displacement jumps across the fracture, respectively, throughout the simulation. The differences between the spring-type and contact-type models are most visible in the tangential direction. Figure~\ref{fig:tangential_all_models} shows that a slightly larger portion of the fracture undergoes tangential deformation at an earlier stage for the spring-type models than for the contact-type models. This behavior is likely related to differences in the tangential constitutive response, as the contact-type models exhibit stick-slip behavior without an intermediate elastic deformation regime. On the other hand, the normal deformation appears unaffected by whether a spring-type or contact-type formulation is used. Instead, the differences in the normal direction are primarily attributed to the choice of elastic normal deformation law (linear or Barton--Bandis). Differences are also observed between the two contact-type models, with larger tangential deformation in the C-Coul-Lin model than in the C-Coul-BB model. This corresponds with the slip tendency reported in Figure~\ref{fig:slip_tendency_model_comparison}, which shows wider stick-regions for the C-Coul-BB model.

\FloatBarrier
\section{Wave Propagation in Rock with Multiple Intersecting Fractures}\label{sec:simulation_example_multiple_intersecting_fractures}
In this section, we investigate the coupled interaction between fracture networks, fracture deformation and elastic wave propagation using the C-Coul-BB fracture deformation model. Unlike the previous examples, which focused on individual fractures or simplified configurations, the present examples consider multiple interacting fractures with different deformation parameters. The purpose of this numerical experiment is to explore the impact of this parameter heterogeneity on both fracture deformation and wave propagation. Such differences may be relevant when interpreting wave propagation and deformation in numerical models of, for example, non-destructive testing, where fractures may exhibit different deformation properties.

We first consider wave propagation in a three-dimensional domain containing both intersecting and isolated fractures to investigate wave–fracture interactions at the scale of the fracture network. A corresponding two-dimensional example is then considered to provide a more detailed fracture-scale view of the deformation response along individual fractures.

\subsection{Three-Dimensional Wave-Fracture Interaction at Network Scale}\label{subsec:three_dimensional_example}
\begin{figure}[pos=ht]
    \centering
    \includegraphics[width=0.5\textwidth]{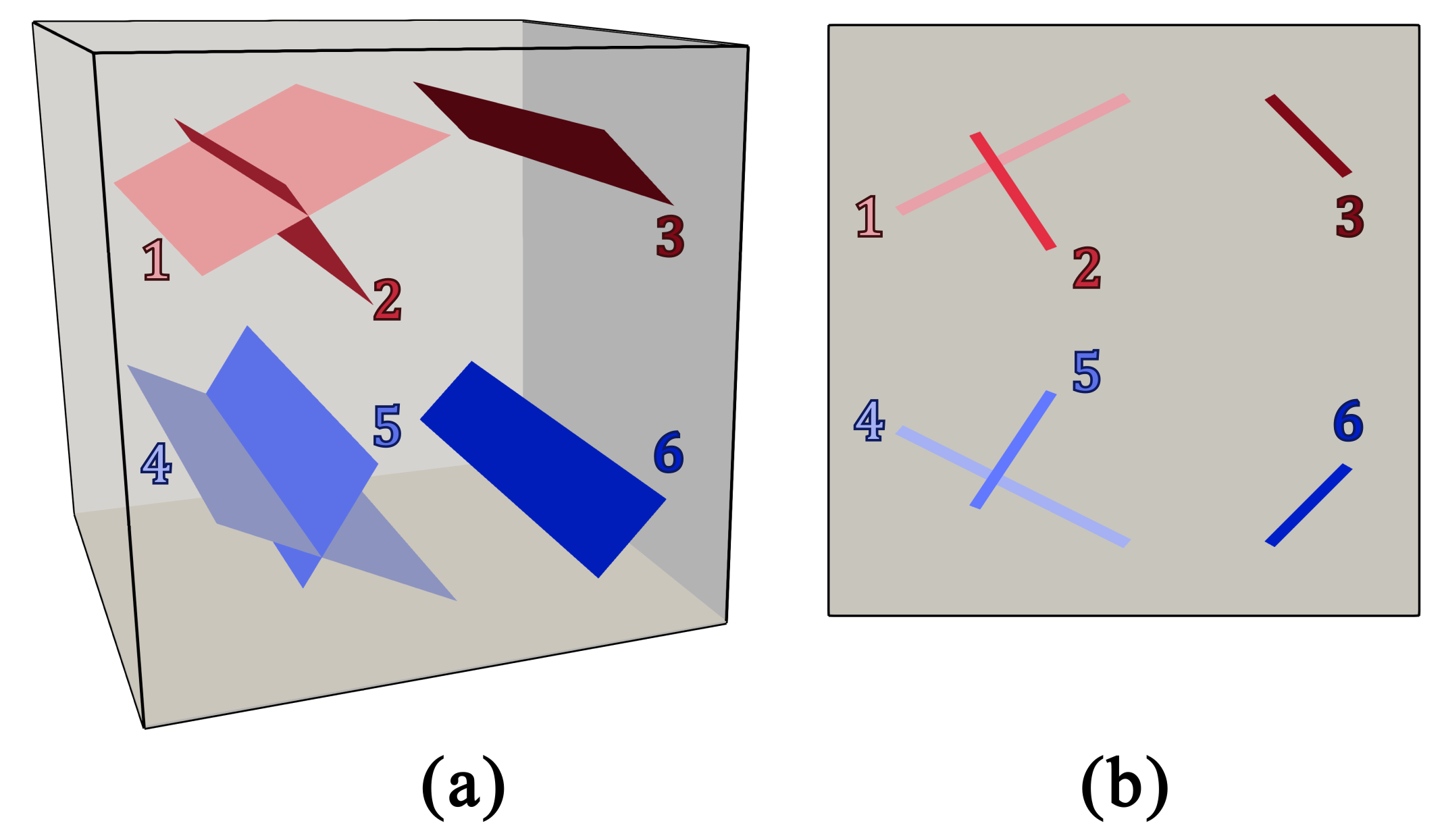}
    \caption{Fracture numbering for the (a) three-dimensional and (b) two-dimensional simulation domains in Section~\ref{sec:simulation_example_multiple_intersecting_fractures}}
    \label{fig:fracture_numbering_3D_2D}
\end{figure}

The three-dimensional domain is a cube with side lengths of $25\,\mathrm{mm}$ and contains six fractures, numbered 1-6, as shown in Figure~\ref{fig:fracture_numbering_3D_2D}a. The domain is obtained by extruding the corresponding two-dimensional geometry shown in Figure~\ref{fig:fracture_numbering_3D_2D}b).
By design, the fracture geometry is symmetric with respect to the $xz$-plane at $y=12.5\,\mathrm{mm}$.

We assign material parameters corresponding to the same soft rock and fractures as considered earlier, namely $\rho = 2600\,\mathrm{kg\,m^{-3}}$, $\mu=\lambda=4.0\,\mathrm{GPa}$ and $K_n=K_\tau=2.0\cdot10^{11}\,\mathrm{Pa\,m^{-1}}$, but with different values of $\Delta u_{\max}$. Fractures 1-3 are assigned $\Delta u_{\max} = 2.0 \cdot 10^{-4}\,\mathrm{m}$, whereas fractures 4-6 are assigned $\Delta u_{\max} = 1.0 \cdot 10^{-5}\,\mathrm{m}$. This results in different effective normal fracture stiffnesses for the upper and lower fractures. Consequently, differences in wave transmissivity and fracture deformation are expected between the two fracture sets.

All boundaries are stress free except the boundary at $x=0$ (left side), which is assigned a time-dependent Dirichlet condition:
\begin{align}
    u(0, y, z,t) =
    \begin{cases}
        \begin{bmatrix}
        U\, w_y(y)\, w_z(z) \sin^2(2\pi f t),\; 0,\; 0
        \end{bmatrix}^{T}
        & \text{for the active region}, \\
        \begin{bmatrix}
        0,\; 0,\; 0
        \end{bmatrix}^{T}
    & \text{otherwise},
    \end{cases}
\end{align}
where $w_y(y)=\sin^2(\pi\frac{y-y_0}{y_1-y_0})$, $w_z(z)=\sin^2(\pi\frac{z-z_0}{z_1-z_0})$ and the active region is a square bounded by $y_0=z_0=6.25\,\mathrm{mm}$ and $y_1=z_1=18.75\,\mathrm{mm}$. The wave has a frequency of $f=100.0\,\mathrm{kHz}$.

We emphasize that the simulation domain is generated by extruding a two-dimensional geometry in the $z$-direction. Consequently, the fracture geometry is invariant along $z$, i.e., the fractures have the same cross-section throughout the extrusion, making the geometry itself quasi-2D. However, the wave is applied only over a subsection of the $x=0$ boundary through the spatial weighting functions $w_y$ and $w_z$, introducing variation in both the $y$- and $z$-directions. As a result, both the wave field and the fracture response vary in all three spatial dimensions, and the simulation therefore represents a fully three-dimensional problem despite the quasi-2D geometry.

\begin{figure}[pos=htbp]
    \centering
    \includegraphics[width=0.5\textwidth]{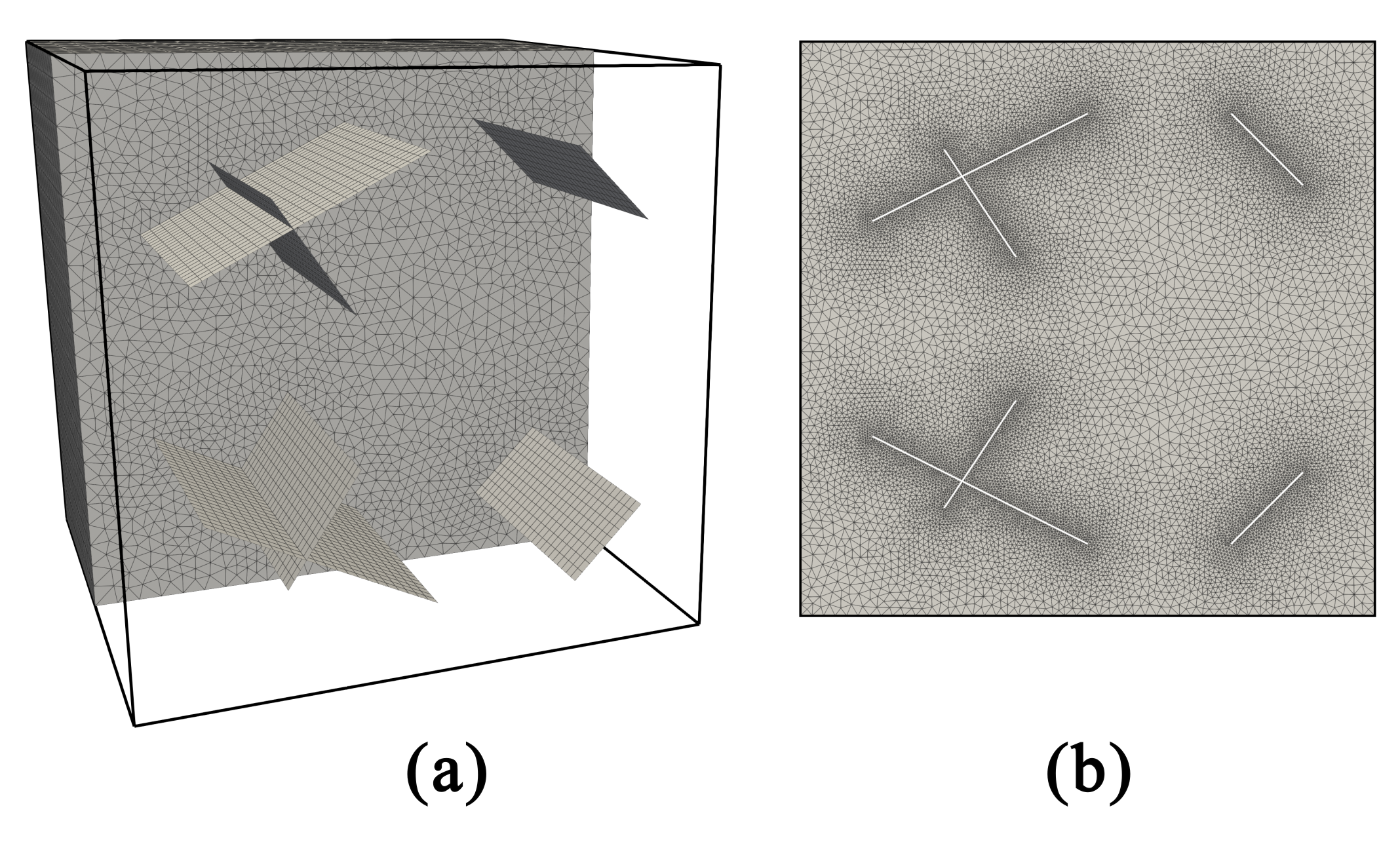}
    \caption{Discretization grids size for the (a) three-dimensional and (b) two-dimensional simulation domains in Section~\ref{sec:simulation_example_multiple_intersecting_fractures}}
    \label{fig:grid_size_3D_2D}
\end{figure}

The spatial and temporal resolutions of the simulation are defined by the mesh size $\Delta x$ and time step $\Delta t$, respectively. The simulation is performed until a final time of $t_{\mathrm{final}} = 13.0\cdot10^{-6}\,\mathrm{s}$ using a time step of $\Delta t = 8.125\cdot10^{-8}\,\mathrm{s}$. Owing to the invariance of the fracture network along the extrusion direction, the three-dimensional domain can be discretized using triangular prism elements, resulting in a lower computational cost than a tetrahedral discretization. The computational mesh has a spatial resolution of $\Delta x = 0.5\,\mathrm{mm}$ in the vicinity of the fractures and $\Delta x = 1.0\,\mathrm{mm}$ near the domain boundaries, as shown in Figure~\ref{fig:grid_size_3D_2D}a. As no convergence analysis was performed for the three-dimensional simulations, the chosen spatial and temporal resolutions are based solely on a criterion of at least 10 points per wave period and wavelength, respectively, to ensure adequate wave resolution.

\begin{figure}[pos=htbp]
    \centering
    \includegraphics[width=0.85\textwidth]{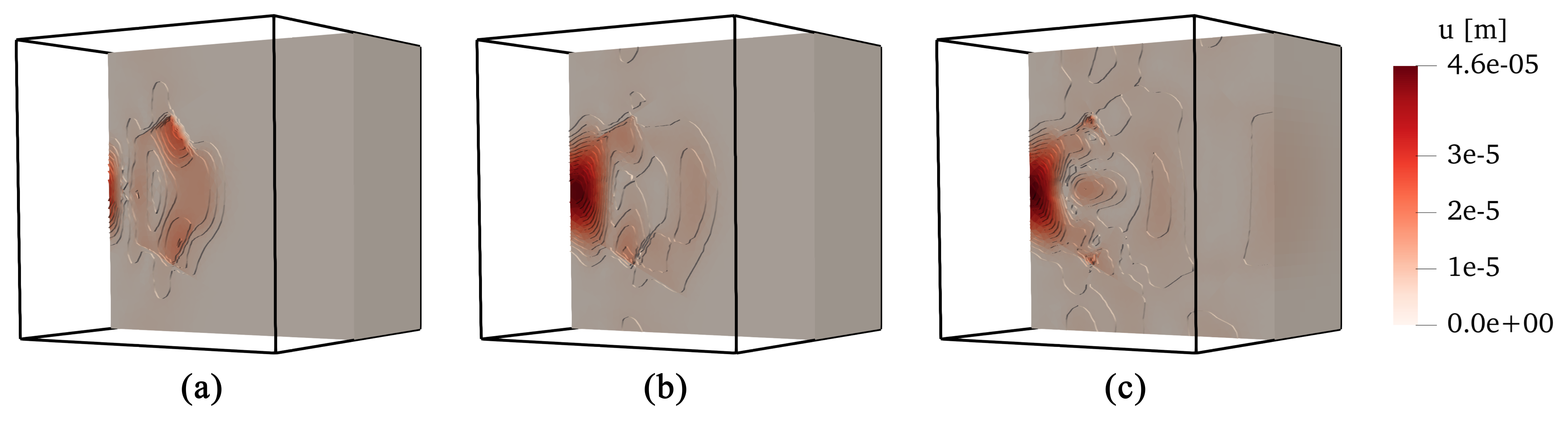}
    \caption{Magnitude of the displacement field at three times during the three-dimensional simulation: (a) $t=6.5\,\mu\mathrm{s}$, (b) $t=8.125\,\mu\mathrm{s}$ and (c) $t=13.0\,\mu\mathrm{s}$. Contours are added for visualization purposes}
    \label{fig:simulation_example_3D_wavefield}
\end{figure}

\begin{figure}[pos=htbp]
    \centering
    \includegraphics[width=0.85\textwidth]{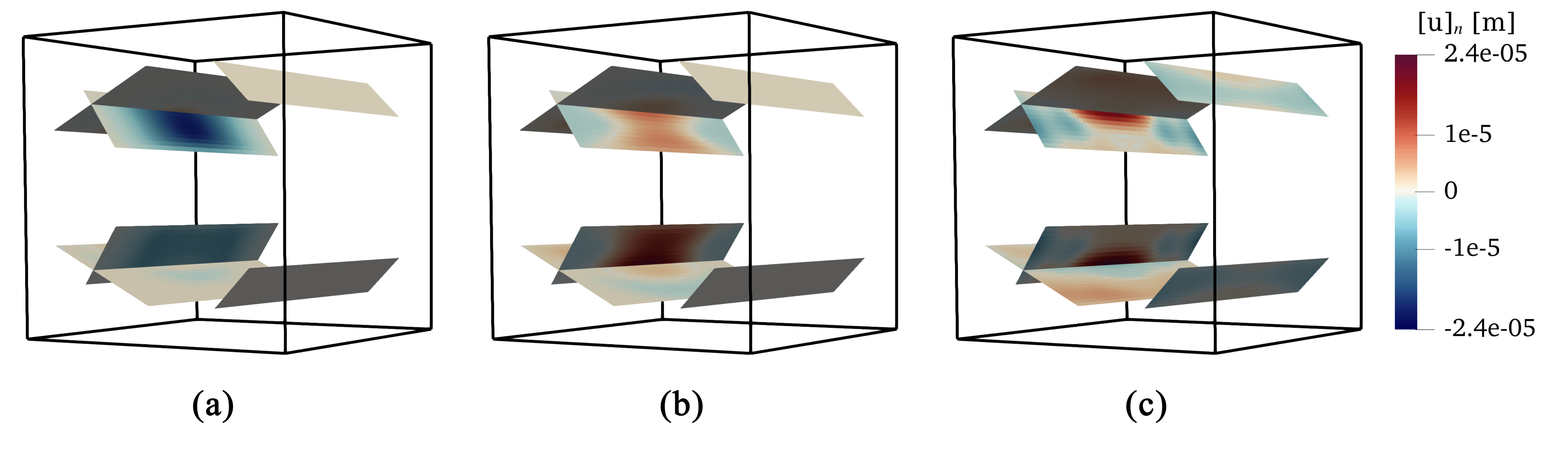}
    \caption{Normal displacement jump $\boldsymbol{\llbracket u \rrbracket}_n$ at three times during the simulation: (a) $t=6.5\,\mu\mathrm{s}$, (b) $t=8.125\,\mu\mathrm{s}$ and (c) $t=13.0\,\mu\mathrm{s}$}
    \label{fig:ex_3d_displacement}
\end{figure}

\begin{figure}[pos=htbp]
    \centering
    \includegraphics[width=0.85\textwidth]{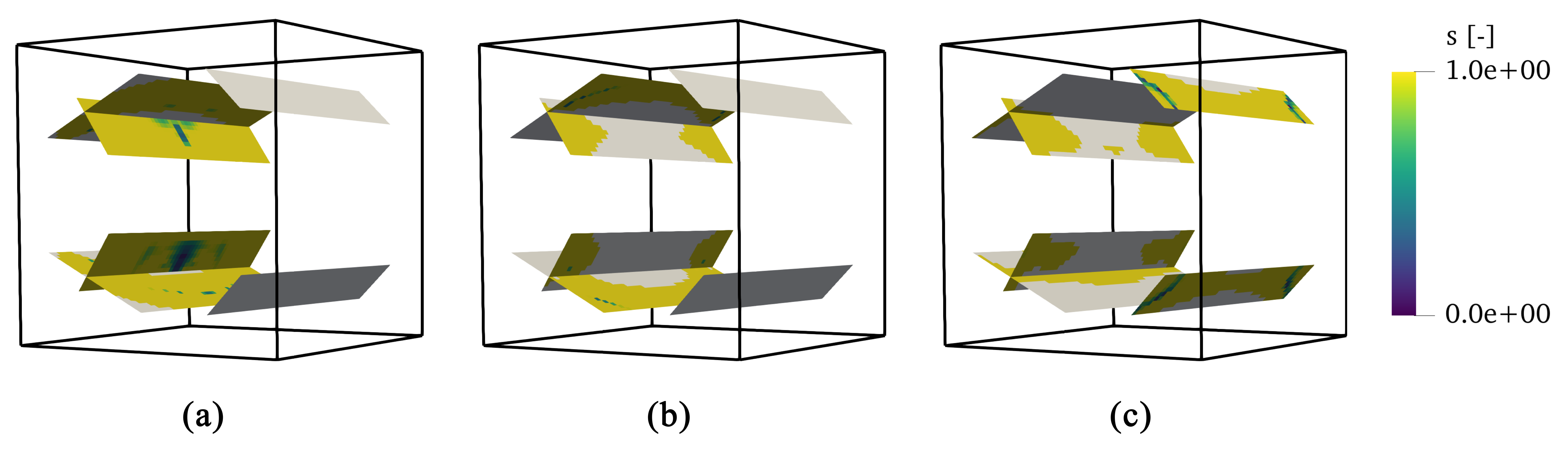}
    \caption{Slip tendency $s$ at three times during the simulation: (a) $t=6.5\,\mu\mathrm{s}$, (b) $t=8.125\,\mu\mathrm{s}$ and (c) $t=13.0\,\mu\mathrm{s}$}
    \label{fig:ex_3d_slip_tendency}
\end{figure}

Simulation results for the wave field, normal displacement jump and slip tendency at three different times are shown in Figures~\ref{fig:simulation_example_3D_wavefield}, \ref{fig:ex_3d_displacement} and \ref{fig:ex_3d_slip_tendency}. Figure~\ref{fig:simulation_example_3D_wavefield} shows that the wave field is affected differently by fractures 1-3 and fractures 4-6 due to the different values of $\Delta u_{\max}$. Specifically, Figure~\ref{fig:simulation_example_3D_wavefield}a) highlights that Fracture 2, relative to Fracture 5, has a lower transmission efficiency and exhibits a higher reflection response. These differences illustrate the influence fracture deformation parameters can have on the wave field in the bulk rock. Further differences between the fractures are observed in Figures~\ref{fig:ex_3d_displacement} and \ref{fig:ex_3d_slip_tendency}, which show distinct normal displacement jumps and slip tendencies for the upper and lower fractures. Additionally, we emphasize the fully three-dimensional nature of the simulation, which is evident from the spatial variation of fracture-related quantities in all coordinate directions across the fracture network.

We also investigate fracture closure for all fractures in the domain. Since compressive fracture states correspond to negative values of the normal displacement jump, the maximum closure of fracture $i$ is quantified as
\begin{equation}
    {\Delta u}_i = \min_{x,y,t}(\boldsymbol{\llbracket u \rrbracket}_{n,i}(x,y,t)), \quad \text{for } i=1,\dots,6,
\end{equation}
where $x$ and $y$ denote coordinates along the fracture surface. The quantity ${\Delta u}_i$ represents the most compressive state attained by fracture $i$ during the simulation and is used to compare the fracture closure of the different fractures relative to $\Delta u_{\max}$.

Table~\ref{tab:fracture_closure_3d} reports the fracture closure values normalized by $\Delta u_{\max}$. The results show that fractures 4-6 are closer to their maximum allowed closure than fractures 1-3, which is expected due to the smaller value of $\Delta u_{\max}$ assigned to fractures 4-6.

The three-dimensional results demonstrate that the two fracture sets (1–3 and 4–6) exhibit distinct deformation and slip characteristics, which in turn give rise to different wave propagation patterns. To complement this network-scale perspective, we next examine an analogous two-dimensional example, which provides a more detailed view of the deformation processes occurring along individual fractures.

\begin{table}[pos=ht]
    \caption{Fracture closure values and closure normalized by the maximum allowed closure}
    \centering
    \begin{tabular}{ccc}
    \hline
    Fracture & ${\Delta u}_i$ & $\lvert {\Delta u}_i \rvert/\Delta u_{\max}$ \\
    \hline
    1 & $-7.277206\times10^{-6}$ & $0.0360$ \\
    2 & $-2.276924\times10^{-5}$ & $0.1138$ \\
    3 & $-5.286359\times10^{-6}$ & $0.0264$ \\
    4 & $-6.742764\times10^{-6}$ & $0.6743$ \\
    5 & $-9.676677\times10^{-6}$ & $0.9677$ \\
    6 & $-6.525257\times10^{-6}$ & $0.6525$ \\
    \hline
    \end{tabular}
    \label{tab:fracture_closure_3d}
\end{table}

\FloatBarrier
\subsection{Two-Dimensional Fracture Response Analysis}\label{subsec:two_dimensional_example}
In order to explore the details of deformation on individual fractures, we now consider a two-dimensional simulation example analogous to that in Section~\ref{subsec:three_dimensional_example}. We refer to Figure~\ref{fig:fracture_numbering_3D_2D}b) for an illustration of the two-dimensional domain. Due to the higher level of spatial refinement attainable in two dimensions, this setup enables a detailed analysis of the evolution of fracture deformation for different $\Delta u_{\max}$ values and fracture orientations. Specifically, we examine the evolution of normal fracture deformation and slip tendency for all fractures throughout the simulation, capturing their variation in space and time.

The boundary conditions are like those shown in Figure~\ref{fig:bc_type_self_convergence} and Equation~\ref{eq:dirichlet_tapered_condition}. The simulation is performed with a grid size of $\Delta x=0.1\,\mathrm{mm}$ near the fractures and $\Delta x=0.5\,\mathrm{mm}$ near the domain boundary (see Figure~\ref{fig:grid_size_3D_2D}b). As in the three-dimensional example, the simulation is run until $t_{final} = 13.0\,\mu\mathrm{s}$ using a time step of $\Delta t = 8.125\cdot10^{-8}\,\mathrm{s}$. Like in Section~\ref{subsec:three_dimensional_example}, both the time step and grid size provide at least 10 points per wave period and wavelength, respectively, to ensure that the wave is resolved. The time step and grid size are also chosen to be within the convergent regime presented in Section~\ref{subsec:fracture_contact_mechanics_friction}.
\begin{figure}[pos=htbp]
    \centering
    \includegraphics[width=0.75\textwidth]{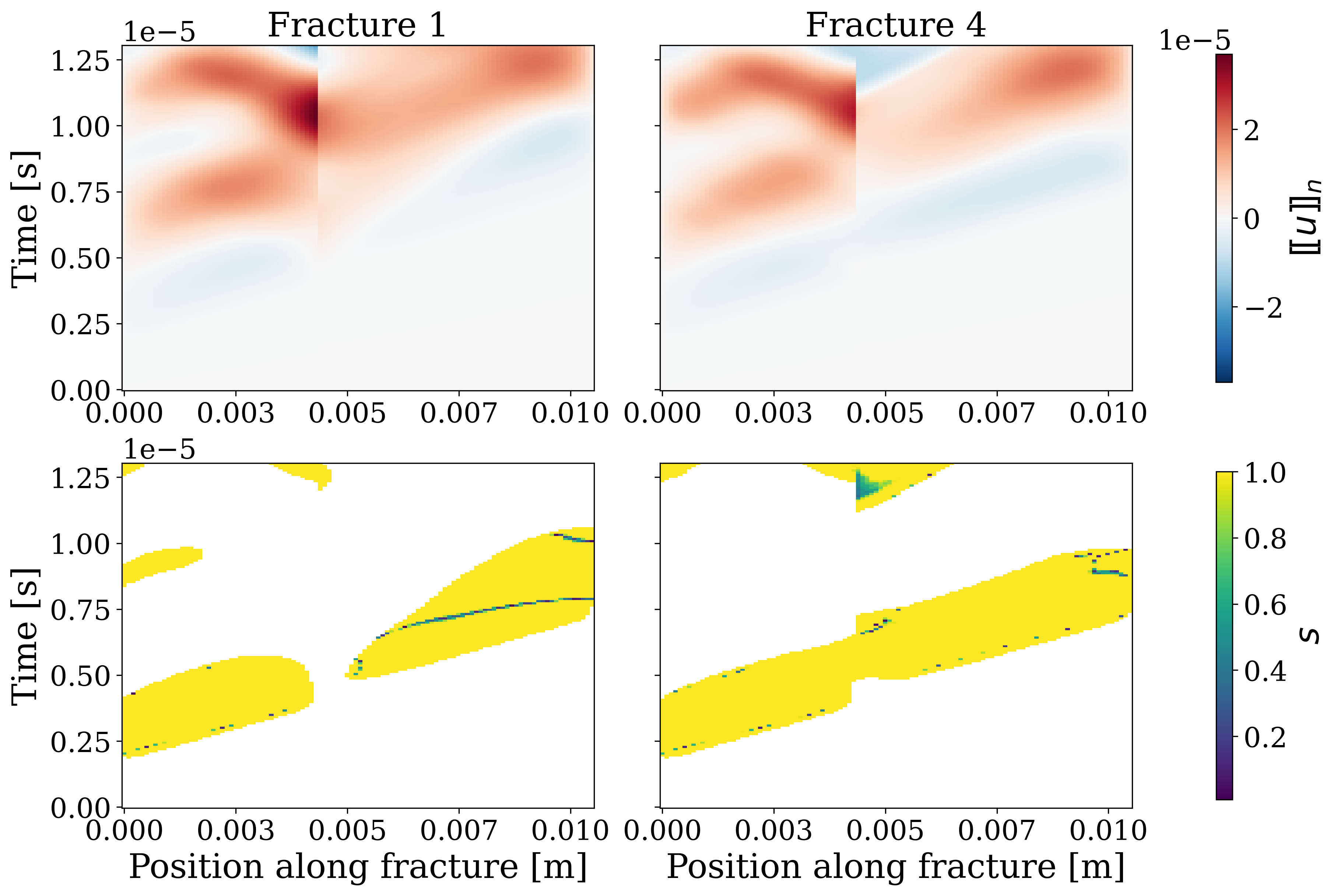}
    \caption{Normal displacement jump $\boldsymbol{\llbracket u \rrbracket}_n$ and slip tendency $s$ for Fractures 1 and 4 in the two-dimensional simulation using the C-Coul-BB model with $\Delta u_{\max}=2.0\cdot10^{-4}\,\mathrm{m}$ (Fracture 1) and $\Delta u_{\max}=1.0\cdot10^{-5}\,\mathrm{m}$ (Fracture 4). White represents $s=\mathrm{NaN}$}
    \label{fig:fractures_1_4}
\end{figure}

\begin{figure}[pos=htbp]
    \centering
    \includegraphics[width=0.75\textwidth]{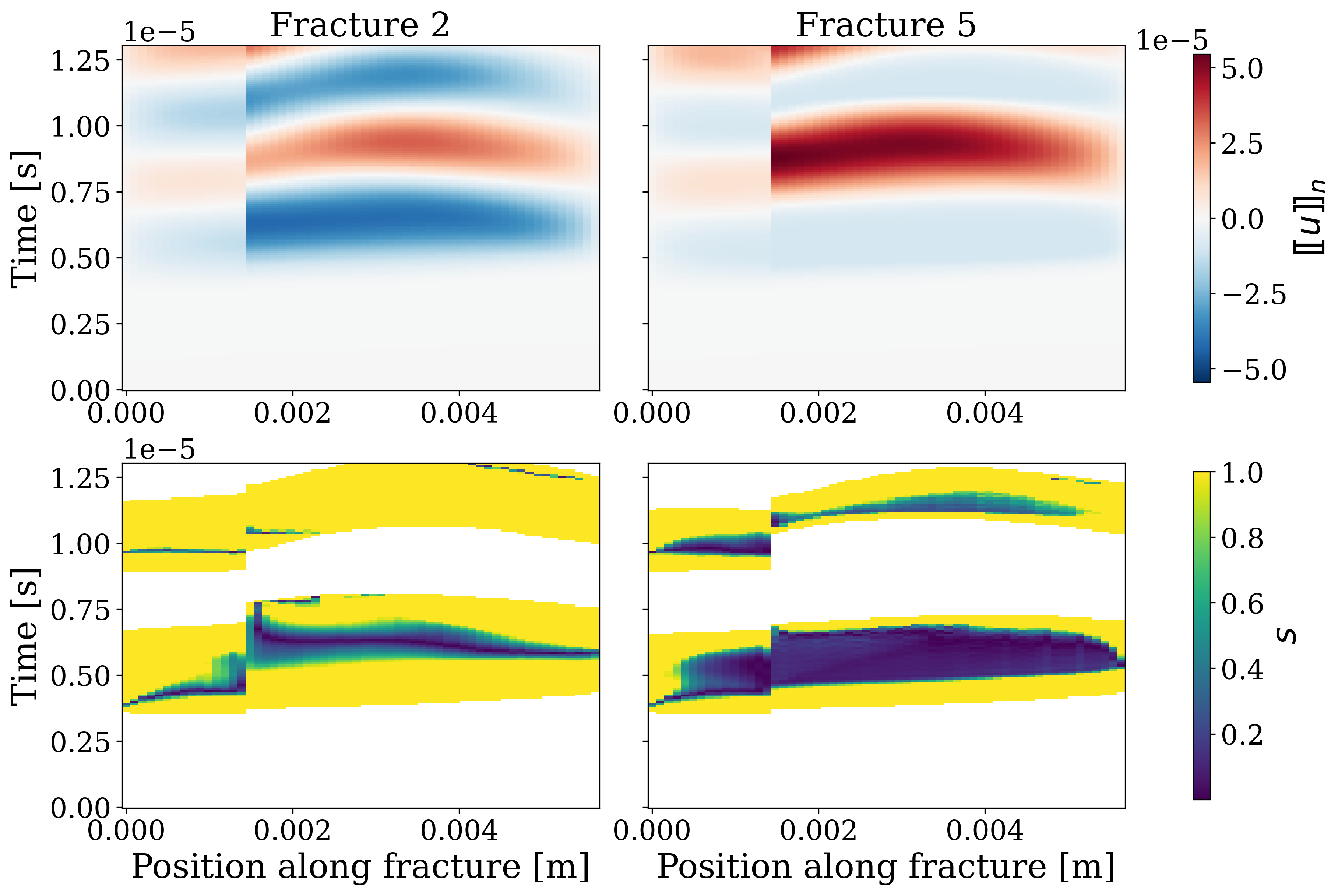}
    \caption{Normal displacement jump $\boldsymbol{\llbracket u \rrbracket}_n$ and slip tendency $s$ for Fractures 2 and 5 in the two-dimensional simulation using the C-Coul-BB model with $\Delta u_{\max}=2.0\cdot10^{-4}\,\mathrm{m}$ (Fracture 2) and $\Delta u_{\max}=1.0\cdot10^{-5}\,\mathrm{m}$ (Fracture 5). White represents $s=\mathrm{NaN}$}
    \label{fig:fractures_2_5}
\end{figure}
\begin{figure}[pos=htbp]
    \centering
    \includegraphics[width=0.75\textwidth]{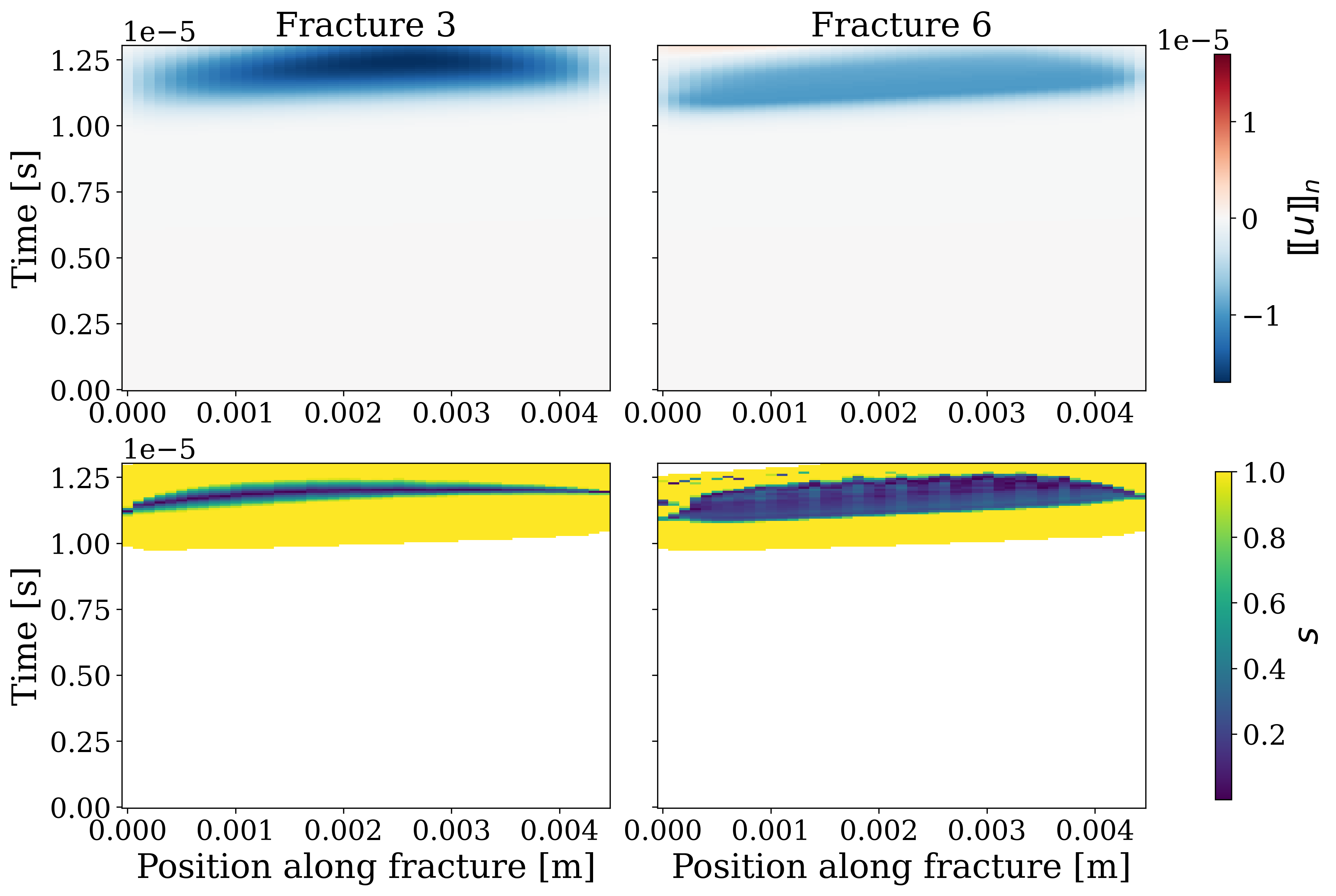}
    \caption{Normal displacement jump $\boldsymbol{\llbracket u \rrbracket}_n$ and slip tendency $s$ for Fractures 3 and 6 in the two-dimensional simulation using the C-Coul-BB model with $\Delta u_{\max}=2.0\cdot10^{-4}\,\mathrm{m}$ (Fracture 3) and $\Delta u_{\max}=1.0\cdot10^{-5}\,\mathrm{m}$ (Fracture 6). White represents $s=\mathrm{NaN}$}
    \label{fig:fractures_3_6}
\end{figure}
\begin{figure}[pos=htbp]
    \centering
    \includegraphics[width=0.85\textwidth]{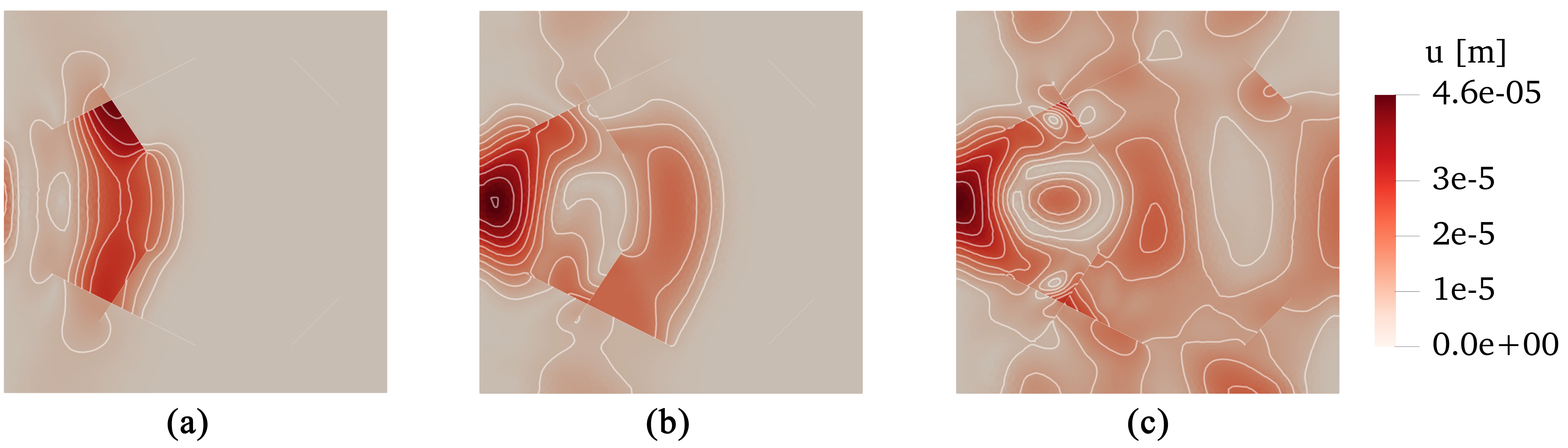}
    \caption{Magnitude of the displacement field at three times during the two-dimensional simulation: (a) $t=6.5\,\mu\mathrm{s}$, (b) $t=8.125\,\mu\mathrm{s}$ and (c) $t=13.0\,\mu\mathrm{s}$. Contour lines are added for visualization purposes}
    \label{fig:wavefield_2d_symmetric_fracture}
\end{figure}

For the fracture deformation results, we refer to Figures~\ref{fig:fractures_1_4}-\ref{fig:fractures_3_6}. Fractures 2 and 5 are oriented such that the incoming wave approaches them at nearly normal incidence. As a result, Fractures 2 and 5 exhibit larger maximum normal deformation than Fractures 1 and 4 which are more closely aligned with the direction of wave propagation. During the initial compression phase ($t \leq 8,\mu\mathrm{s}$), Figure~\ref{fig:fractures_2_5} shows that Fracture 5 exhibits a more pronounced stick behavior than Fracture 2, likely due to the stronger closure associated with its lower $\Delta u_{\max}$. Similarly, Fracture 6 exhibits more stick behavior than Fracture 3 (Figure~\ref{fig:fractures_3_6}). This is consistent with more efficient wave transmission through Fractures 4 and 5 than through Fractures 1 and 2, which results in a larger-amplitude incident wave at Fracture 6 than at Fracture 3. The wave field is shown in Figure~\ref{fig:wavefield_2d_symmetric_fracture}, where an asymmetry similar to that observed in Figure~\ref{fig:simulation_example_3D_wavefield} is seen. This further demonstrates, in a two-dimensional setting, that fracture deformation parameters influence wave propagation in the surrounding rock.

The results in this section show that the model can resolve the evolution of normal fracture deformation and slip tendency in space and time. We also observe that relatively small variations in the parameter $\Delta u_{\max}$ lead to noticeably different deformation responses. Although $\Delta u_{\max}$ governs normal fracture closure, its influence also extends to the tangential fracture response. Furthermore, variations in $\Delta u_{\max}$ also significantly affect the wave field, as observed in both the two- and three-dimensional simulations.

\FloatBarrier
\section{Conclusion}\label{sec:conclusion}
We have presented a finite volume approach, namely the MPSA-Newmark method, for elastic wave propagation in deformable fractured media. The approach accommodates fracture deformation models with increasing physical complexity, ranging from linear and nonlinear spring-type formulations to frictional contact mechanics.

Numerical convergence analyses were performed in space and time for all fracture deformation models, including configurations with diagonal fractures and frictional effects. The observed convergence rates are within the expected range of 1-2, demonstrating the numerical convergence of the proposed methodology across a broad range of fracture behaviors.

A simulation example in a two-dimensional fractured medium was used to illustrate the influence of the fracture deformation model on both wave propagation and fracture response. The results demonstrate that different fracture deformation models can produce markedly different wave transmission, scattering and deformation patterns, emphasizing the importance of fracture constitutive behavior in wave propagation through fractured media.

Finally, two- and three-dimensional simulation examples with multiple intersecting fractures were presented, all governed by the frictional contact model with Barton-Bandis elastic normal deformation. The results show that variations in the Barton--Bandis parameter $\Delta u_{\max}$ can lead to substantial differences in wave transmission, fracture deformation and slip tendency. These findings highlight the importance of fracture constitutive properties for wave propagation in fractured media and demonstrate that nonlinear normal deformation and frictional contact can significantly influence both the wave field and the mechanical response of fractures.

Overall, our experience indicates that the proposed method provides robust approximations of physical quantities of interest in both two- and three-dimensional fractured geometries. The results further highlight the importance of selecting fracture constitutive models and parameters that reflect the physical mechanisms of interest, as these choices can substantially affect wave transmission, scattering and fracture deformation. The proposed methodology thus provides a tool for investigating the coupling between fracture deformation and wave propagation in fractured rocks. By capturing these interactions while exhibiting no apparent spurious oscillations in the wave fields of the presented simulations, the method can support the interpretation of wave patterns in applications such as non-destructive testing. Furthermore, the compatibility of the finite volume approach with subsurface flow models provides a natural foundation for future extensions to poroelastic modeling in fractured media, including elastic wave propagation.

\appendix
\numberwithin{equation}{section} 

\section{Errors in cell and face quantities measured against known analytical solutions}\label{appendix_known_analytical_solution}
All numerical errors reported here are discrete $L^2$-type norms. To present the norms, we first cover the necessary notation for the discrete representation of the domain. The computational domain is divided into the fracture subdomain $\Omega_l$ and the rock matrix subdomain $\Omega_h$. As the concept is the same for $\Omega_l$ and $\Omega_h$, we omit the subscripts $l$ and $h$ and instead use $\Omega$ to keep a cleaner presentation. All subdomains are divided into $k$ non-overlapping cells $K_i$ for $i=1,2,\text{…},k$ such that $\Omega = \bigcup_i \overline{K_i}$. We denote the set of cells in a subdomain by $\mathcal{T}$, a face on the boundary of a cell by $f$ and the set of cell-faces in a subdomain by $\mathcal{F}$. The geometric measure of a particular cell $K$ is denoted by $m_K$, whereas that of a particular face $f$ is denoted by $m_f$. 

We consider now generic cell and face quantities denoted by $\xi$ and $\chi$, respectively, and use boldface to distinguish the discrete solution from the continuous solution. The weighted norms for generic discrete cell and face quantities $\boldsymbol{\xi}$ and $\boldsymbol{\chi}$ are: 
\begin{align}
    \lVert\boldsymbol{\xi}\rVert_\mathcal{T}&=\left(\sum_{K\in \mathcal{T}} m_K (\boldsymbol{\xi}_K\cdot \boldsymbol{\xi}_K)\right)^{1/2}, \label{eq:cell_norm} \\
    \lVert\boldsymbol{\chi}\rVert_\mathcal{F}&=\left(\sum_{f\in \mathcal{F}} \frac{1}{D} m_f (d_L - d_R)\cdot n_f (\boldsymbol{\chi}_f\cdot \boldsymbol{\chi}_f)\right)^{1/2}. \label{eq:face_norm}
\end{align}
The symbols $d_L$ and $d_R$ denote the vectors pointing from the face-center $x_f$ to the cell-center of the two cells $L$ and $R$ on either side of the face $f$. The symbol $n_f$ denotes the face normal vector pointing from cell $R$ to cell $L$, and $D$ represents the dimension of $\Omega$. If $f$ is a face on the boundary or a face coinciding with a fracture cell, one of $d_L$, $d_R$ is the zero vector.

We compute the numerical errors of cell- and face-centered quantities relative to the analytical solution projected onto cell-centers and face-centers. For this we introduce the projection $\Pi_\mathcal{T}$ which returns a cell-centered quantity ($(\Pi_{\mathcal{T}}\xi)_K=\xi(x_K)$) and the projection $\Pi_\mathcal{F}$ which returns a face-centered quantity ($(\Pi_{\mathcal{F}}\chi)_f=\chi(x_f)$). Then, the relative errors reported in this section are computed according to:
\begin{align}
    \mathcal{E}_\xi &= \frac{\lVert \boldsymbol{\xi} - \Pi_\mathcal{T}\xi \rVert_\mathcal{T}}{\lVert \Pi_\mathcal{T}\xi \rVert_\mathcal{T}}, \label{eq:error_cell} \\
    \mathcal{E}_\chi &= \frac{\lVert \boldsymbol{\chi} - \Pi_\mathcal{F}\chi \rVert_\mathcal{F}}{\lVert \Pi_\mathcal{F}\chi \rVert_\mathcal{F}}. \label{eq:error_face}
\end{align}
Subscript $\xi$ and $\chi$ on the error symbol $\mathcal{E}$ denote error in cell-quantity and face-quantity, respectively. All errors are computed at the final time in the simulation.

\section{Error between numerical solutions on non-matching grids}\label{appendix_nonmatching_grids}
To compare solutions on two non-matching grids $\mathcal{T}_h$ (fine) and $\mathcal{T}_H$ (coarse), we first map the cell-centered vector solution $\boldsymbol{\xi}$ to a piecewise linear representation $\xi_{\mathcal{R}}^D$ where $D \in \{\mathrm{1}, \mathrm{2}\}$ denotes the grid dimension and $\mathcal{R}$ indicates that the field is reconstructed. The linear representation is based on nodal values which are defined as a weighted average of cell-centered values:
\begin{align}
    \boldsymbol{\zeta}_N = \frac{1}{\sum_{K \in \mathcal{T}_N} m_K}
    \sum_{K \in \mathcal{T}_N} m_K \boldsymbol{\xi}_K,
\end{align}
where $\mathcal{T}_N$ denotes the set of cells sharing node $N$ and $m_k$ is the geometric measure of cell $K$. The reconstructed field $\xi_{\mathcal{R}}^D$ is obtained by linear interpolation of the nodal values $\boldsymbol{\zeta}_N$ on each cell:
\begin{align}
    \xi_{\mathcal{R},K}^{\mathrm{2}}(x,y) = a_K x + b_K y + c_K,
    \qquad
    \xi_{\mathcal{R},K}^{\mathrm{1}}(s) = d_K s + e_K.
\end{align}
All the coefficients are determined by enforcing that the nodal values are reproduced on each cell vertex. The reconstruction procedure is applied for each component of the vector solution. Combining the linear solution on all cells yields the global piecewise linear representations $\xi_{\mathcal{R}}^1$ and $\xi_{\mathcal{R}}^2$.

Following \citet{Varela2025, Varela2026Zenodo}, we introduce a transfer grid $\mathcal{T}_{Hh}$ obtained by subdividing the non-empty intersections of cells in $\mathcal{T}_h$ and $\mathcal{T}_H$. Each $K_{Hh} \in \mathcal{T}_{Hh}$ is contained in exactly one fine cell $K_h$ and exactly one coarse cell $K_H$. Since the reconstructed fields are piecewise linear, their restriction to each transfer element is also linear. 

The error between the reconstructed fields $\xi_{\mathcal{R}}^{D,h}$ (fine grid) and $\xi_{\mathcal{R}}^{D,H}$ (coarse grid) is defined in a relative $L^2$-norm as
\begin{align}
    E =
    \left(
    \frac{
    \sum_{K_{Hh} \in \mathcal{T}_{Hh}}
    \int_{K_{Hh}} \left\| \xi_{\mathcal{R}}^{D,H} - \xi_{\mathcal{R}}^{D,h} \right\|^2 dx
    }{
    \sum_{K_{Hh} \in \mathcal{T}_{Hh}}
    \int_{K_{Hh}} \left\| \xi_{\mathcal{R}}^{D,h} \right\|^2 dx
    }
    \right)^{1/2} \quad \text{for }D \in \{\mathrm{1}, \mathrm{2}\}.
    \label{eq:l2_error_nonmatching_grids}
\end{align}
Since the integrands in Equation~\ref{eq:l2_error_nonmatching_grids} are quadratic polynomials, each element contribution is evaluated exactly using the analytical finite element mass matrix for linear basis functions.

\printcredits

\section*{Declaration of competing interest}
The authors declare that they have no known competing financial interests or personal relationships that could have appeared to 
influence the work reported in this paper.

\section*{Acknowledgements}
This project has received funding from the European Research Council (ERC) under the European Union’s Horizon 2020 research and innovation program (grant agreement No. 101002507). JMN and BW acknowledge the support of the Centre for Advanced Study in Oslo, Norway that funded and hosted the research project Mathematical Challenges in Brain Mechanics during the academic year of 2025/2026.

\section*{Declaration of generative AI and AI-assisted technologies in the manuscript preparation process}
During the preparation of this work the author(s) used ChatGPT and Microsoft Copilot in order to refine the academic language. After using this tool/service, the author(s) reviewed and edited the content as needed and take(s) full responsibility for the content of the published article.

\section*{Data availability}
The source code for all simulations is open-source and is available at Zenodo in \citet{Jacobsen2026Zenodo}.

\bibliographystyle{cas-model2-names}

\bibliography{cas-refs}

\end{document}